\documentclass{article}

\usepackage[leqno,intlimits]{amsmath}
\usepackage{amssymb}
\usepackage{amsthm}
\usepackage{hyperref}
\usepackage{appendix}
\usepackage{calc}
\usepackage{pstricks}
\usepackage{multido}
\psset{unit=1pt}
\psset{linewidth=0.5pt}

\newcounter{z}

\newcommand*{\hop}{\bigskip\noindent}
\newcommand*{\wih}{\widehat}
\newcommand*{\wt}{\widetilde}
\newcommand*{\vr}{\varrho}
\newcommand*{\la}{\lambda}

\newcommand*{\de}{\delta}
\newcommand*{\te}{\theta}
\newcommand*{\ve}{\varepsilon}
\newcommand*{\vp}{\varphi}
\newcommand*{\ze}{\zeta}
\newcommand*{\uetap}{\underline{\eta}^+}
\newcommand*{\uetae}{\underline{\eta}^\text{eq}}
\newcommand*{\etae}{\eta^\text{eq}}
\newcommand*{\al}{\alpha}
\newcommand*{\ga}{\gamma}
\newcommand*{\Zb}{\mathbb Z}
\newcommand*{\Nb}{\mathbb N}
\newcommand*{\Rb}{\mathbb R}
\newcommand*{\Hc}{\mathcal H}
\newcommand*{\Mc}{\mathcal M}
\newcommand*{\Ic}{\mathcal I}
\newcommand*{\Jc}{\mathcal J}
\newcommand*{\un}[1]{\underline{#1}}
\newcommand*{\om}{\omega}
\newcommand*{\omp}{\omega}
\newcommand*{\omm}{\omega^-}
\newcommand*{\ome}{\omega^\text{eq}}
\newcommand*{\uomp}{\underline\omega}
\newcommand*{\uomm}{\underline{\omega}^-}
\newcommand*{\uome}{\underline{\omega}^\text{eq}}
\newcommand*{\omin}{\om^{\text{min}}}
\newcommand*{\omax}{\om^{\text{max}}}
\newcommand*{\be}{\begin{equation}}
\newcommand*{\ee}{\end{equation}}
\newcommand*{\ba}{\begin{aligned}}
\newcommand*{\ea}{\end{aligned}}
\newcommand*{\Ev}{{\bf E}}
\newcommand*{\Pv}{{\bf P}}
\newcommand*{\Vv}{{\text{\bf Var}}}
\newcommand*{\Cov}{{\text{\bf Cov}}}
\newcommand*{\fl}[1]{\lfloor{#1}\rfloor}
\newcommand*{\e}[1]{\text{\rm e}^{#1}}
\newcommand*{\di}{\,\text{\rm d}}
\newcommand*{\dii}{\text{d}}

\newtheorem{tm}{Theorem}[section]
\newtheorem{lm}[tm]{Lemma}

\newtheorem{cor}[tm]{Corollary}
\newtheorem{pr}[tm]{Proposition}
\newtheorem{ass}[tm]{Assumption}

\providecommand{\abs}[1]{\lvert#1\rvert}
\def\ind{\mathbf{1}}
\def\Zca{Y} %%normalization for an exponential integral in an appendix

\numberwithin{equation}{section}

\begin{document}
\title{Microscopic concavity and fluctuation bounds in a class of deposition processes}
\author{
M\'arton Bal\'azs\thanks{Budapest University of Technology and Economics. Part of this work was done while M.\ Bal\'azs was affiliated with the MTA-BME Stochastics Research Group},
J\'ulia Komj\'athy\thanks{Budapest University of Technology and
Economics},
Timo Sepp\"al\"ainen\thanks{University of Wisconsin-Madison\newline
M.\ Bal\'azs and J. Komj\'athy were partially supported by the Hungarian
Scientific Research Fund (OTKA) grants K60708, TS49835, F67729, and by Morgan Stanley Mathematical Modeling Center.
M.\ Bal\'azs was also partially funded by the Bolyai Scholarship of the Hungarian Academy of Sciences.
T.\ Sepp\"al\"ainen was partially supported by National Science Foundation grant DMS-0701091 and by the Wisconsin Alumni Research Foundation.}}
\maketitle
\begin{abstract}
We prove fluctuation bounds for the particle current in totally asymmetric zero range processes in one dimension with nondecreasing, concave jump rates whose slope decays exponentially. Fluctuations in the characteristic directions have order of magnitude $t^{1/3}$. This is in agreement with the expectation that these systems lie in the same KPZ universality class as the asymmetric simple exclusion process. The result is via a robust argument formulated for a broad class of deposition-type processes. Besides this class of zero range processes, hypotheses of this argument have also been verified in the authors' earlier papers for the asymmetric simple exclusion and the constant rate zero range processes, and are currently under development for a bricklayers process with exponentially increasing jump rates.
\end{abstract}

\noindent {\bf Keywords:} Interacting particle systems, universal fluctuation bounds, $t^{1/3}$-scaling, second class particle, convexity, asymmetric simple exclusion, zero range process

\hop
{\bf 2000 Mathematics Subject Classification:} 60K35, 82C22

\section{Introduction}

This paper studies anomalous current fluctuations of attractive interacting systems in one dimension with one conserved quantity. The family of models considered includes the asymmetric exclusion, the zero range, misanthrope-type and many other processes. In the asymmetric case (to be specified later) the Eulerian scaling of such a system leads to a (deterministic) hyperbolic conservation law with a hydrodynamic flux function \(\Hc(\vr)\). The \emph{characteristics} of the conservation law is of particular importance both for the PDE itself and for the underlying stochastic system. Recently, the current fluctuations through the characteristic lines drew much attention. The behavior of these fluctuations is fundamentally determined by the form of \(\Hc\). Rigorous results exist for examples that fall in two categories.
 
 \medskip
 
\emph{Order $t^{1/4}$ fluctuations.}
 When $\Hc$ is linear the fluctuations are of order $t^{1/4}$ and converge to
Gaussian processes related to fractional Brownian motion.  This has been proved for 
independent particles \cite{dgl, kumarrw, flucha} and the random average process 
\cite{raprwre, rap}.  
 
 \medskip
 
\emph{Order $t^{1/3}$ fluctuations.}
 When $\Hc''(\vr)\ne 0$ the fluctuations are of order $t^{1/3}$ and converge to
distributions and processes related to the Tracy-Widom distributions from random
matrix theory.  The most-studied examples are the totally asymmetric
simple exclusion process (TASEP), the polynuclear growth model
and the Hammersley process.  Two types  of mathematical work should be distinguished.

(a) Exact limit
distributions have been derived with techniques of asymptotic analysis applied to determinantal
representations of the probabilities of interest.
Most of this work has dealt with particular deterministic initial conditions, and 
 the stationary situation has been less studied. 
The seminal results appeared in \cite{bdj} for the last-passage 
version of the Hammersley process and in \cite{1/3} for the last-passage model
associated with TASEP.  Current fluctuations for stationary TASEP were 
analyzed in \cite{ferspohn}.  Here is a selection of further results in this 
direction: \cite{baik-rain-00, boro-etal-07, grav-trac-wido-01, joha-03, prah-spoh-02}. Recently, the asymmetric simple exclusion (ASEP) also got within reach of these techniques \cite{tr-w-curr-asep}.

(b) Probabilistic approaches exist to prove fluctuation bounds
of the correct order.  The seminal work  \cite{cuberoot} was on 
  the last-passage version of the Hammersley process, and then the
  approach was adapted to the last-passage model associated with TASEP  \cite{third}.
  The next step was the development of a proof that works for particle systems:
   the ASEP was treated  in \cite{se2/3} and 
the totally asymmetric zero range process (TAZRP) with constant jump rate in \cite{julizrp}.
The ASEP work \cite{se2/3} was the first to prove $t^{1/3}$ order of 
  fluctuations for a process  where particle motion is
not restricted to totally asymmetric. Resolvent methods were also applied in \cite{quava2,quava} to extend the results from nearest neighbor ASEP to exclusion processes with non nearest neighbor jumps.
  
 \medskip

The present paper takes a further step toward universality of the $t^{1/3}$ order of fluctuations in the case $\Hc''(\vr)\ne 0$. We rewrite our earlier proof for ASEP and constant rate TAZRP in a fairly general way, extract and formulate in general terms a particular feature of these two models that made our proof work. For reasons to be explained later we call this feature \emph{microscopic concavity}. With this notion in hand we extend the \(t^{1/3}\) scaling result for a class of totally asymmetric zero range processes (with non constant rates). We remark at this point that jump rates of this example have a much richer behavior than the constant rates of those featured in anomalous scaling proofs so far. Further generalizations now only require the verification of microscopic concavity. Product form invariant distributions are critically important for the method.

The hypothesis of microscopic concavity consists of control of second class particles that is a microscopic counterpart of the macroscopic effect that concavity of $\Hc$ has on characteristics. We make this technically precise in Section \ref{sc:microconc}. Once the microscopic concavity assumption is made the proof works for the entire class of processes. This then is the sense in which we take a step toward universality. As a by-product, we also obtain superdiffusivity of the second class particle in the stationary process.

Earlier proofs of $t^{1/3}$ fluctuations have been quite rigid
 in the sense that they work only for particular 
 cases of the models where special combinatorial properties emerge as if 
 through some   fortuitous coincidences. There is basically no room for perturbing
 the rules of the process.   By contrast, the proof given in 
 the present paper works for the whole class of processes.  The hypothesis of
 microscopic concavity that is required is certainly nontrivial.   But it does not seem
 to  rigidly exclude all  but a handful of the processes in the 
 broad class.  
 The estimates that it requires can probably be proved in different ways for
 different subclasses of the processes. And the proof itself may evolve further
and weaken the  hypothesis required. 

\medskip

To summarize, we are currently
 able to verify the required  hypothesis of microscopic concavity  for  the 
 following three subclasses of processes.
 
(i) The asymmetric simple exclusion process (ASEP).  Full details of this case are reported
elsewhere \cite{asepsimple} and we give a brief informal description in Section \ref{sc:asep}.
This proof is somewhat simpler than the earlier one given in \cite{se2/3}. 

(ii) Totally asymmetric zero range processes (TAZRP) with a concave jump rate function whose slope 
decreases geometrically, and may be eventually constant.  This example has been out of reach for existing methods, so it is completely new in this context. It is developed
fully in the present paper. As a special case, the result of \cite{julizrp} for the constant rate TAZRP is also recovered.

(iii) The totally asymmetric bricklayers process with convex, exponential jump rate.
This system satisfies the analogous \emph{microscopic convexity}. 
Due to the fast growth of the jump rate function this example needs more
preliminary work than was sensible to include in the present paper, and so the
result will be published separately in the future.  

We expect that a broader class of totally asymmetric 
concave zero range processes should be amenable
to further progress because a key part of the hypothesis can be verified, and only
a certain tail estimate is missing.  We explain this in Section \ref{sc:1zrp}. 

Interacting particle systems can naturally be given a surface growth representation where integrated particle current becomes the height of a surface and particle occupations become (negative) discrete gradients of this surface. We found this picture extremely helpful in visualizing currents and couplings, hence this is the way we introduce and handle the processes.

\medskip

This paper has two parts.  In the main part we prove the general fluctuation bound under
the assumptions needed for membership in the class of processes and 
the assumption of microscopic concavity.  
 The remainder of the paper  shows that  the
assumptions required by the general result are satisfied by a class of zero range
processes.  Here is a section by section outline. 

 In  Section \ref{sc:modres}
  we define the general family of processes under consideration, describe the 
  microscopic concavity property and other assumptions used, and state the
  general results.  Partly as corollaries to the fluctuation bound along the
  characteristic we obtain a law of large numbers for a second class particle 
  and limits that show how fluctuations in non-characteristic directions on the 
  diffusive scale   come directly from fluctuations of the initial state.
     Section \ref{sc:ex} describes two examples.
   Section \ref{sc:asep} gives a brief description of how    the asymmetric simple exclusion process 
     (ASEP)  satisfies the assumptions of our general theorem.  (Full details for this
     example are reported in \cite{asepsimple}.)   Section \ref{sc:1zrp} describes
     a class of totally asymmetric zero range processes with  concave jump rates that increase with exponentially decaying slope.
     
     The general theorem is proved   in two parts: the upper bound in Section \ref{sc:ub} and 
     the lower bound in Section \ref{sc:lb}. Section \ref{sc:Qlln} proves 
a strong law for the second class particle, partly as a corollary of the main fluctuation bounds.     
      We then return to the zero range example
     and give a complete proof for this class of processes in  Section \ref{sc:zrp}.   
     
      A three-part Appendix contains auxiliary computations for the stationary distribution
     and   hydrodynamic flux function.  In particular,   if the jump rate function of a zero 
     range process is concave and not
     linear then the hydrodynamic flux $\Hc$ satisfies $\Hc''(\vr)<0$ for all 
     densities $0<\vr<\infty$.  

\subsubsection*{Notation} We summarize here some notation for easy reference.
$\Zb^+=\{0,\,1,\,2,\dotsc\}$, $\Rb^+=[0,\,\infty)$. Centering a random variable is
denoted by $\wt X=X-EX$.  Constants \(C_{\centerdot},\,\al_{\centerdot}\) 
do not depend on time, but may depend on the density parameter \(\vr\) and their values can change from line to line.
The numbering of these constants is of no particular significance and is meant only to facilitate
following the arguments.  

\section{Definitions and results}\label{sc:modres}

We  define the class of processes studied in this paper, give a list of
examples, and   discuss
some of basic properties.  Then come the hypotheses and main
results of this paper, followed by two examples of subclasses of processes 
for which the hypotheses can be verified.

\subsection{A family of deposition processes}\label{sc:modfam}

The family of processes we consider is the one described in \cite{varj2nd}, 
and we repeat the definition here. We start with the interface growth picture, 
but we end up using the height and particle languages interchangeably.
  For extended-integer-valued boundaries 
 $-\infty\le\omin\le0$ and $1\le\omax\le\infty$ 
define the single-site state space
\[
I:\,=\left\{z\in\Zb\,:\,\omin-1<z<\omax+1\right\}
\]
and the increment configuration space
\[
\Omega:\,=\left\{\un\om=(\om_i)_{i\in\Zb}\ :\ \om_i\in I\right\}=I^{\Zb}.
\]
At times it will be convenient to have notation for the increment configuration
$\un\de_i\in\Omega$ with exactly one nonzero entry equal to $1$:
\be
(\un\de_i)_j=\left\{\ba
&1,&&\text{for }i=j,\\
&0,&&\text{for }i\ne j.
\ea\right.\label{eq:dedef}
\ee

For each pair of neighboring sites $i$ and $i+1$ of $\Zb$  imagine a column of bricks over
the interval $(i,\,i+1)$. The height $h_i$ of this column is integer-valued.  The components of a  
  configuration $\un{\om}\in\Omega$ are  the negative discrete gradients of the heights: 
  $\om_i=h_{i-1}-h_i\,\in I$.
  
   The evolution is described by jump processes whose rates $p$ and $q$ are nonnegative 
   functions on $I\times I$. 
Two types of moves are possible. A brick can be deposited:
\be
\left.
\ba
\left(\om_i,\,\om_{i+1}\right)&\longrightarrow\left(\om_i-1,\,\om_{i+1}+1\right)\\
h_i&\longrightarrow h_i+1
\ea
\right\}
\text{with rate}\ p(\om_i,\,\om_{i+1}),\label{eq:add}
\ee
or removed:
\be
\left.
\ba
\left(\om_i,\,\om_{i+1}\right)&\longrightarrow\left(\om_i+1,\,\om_{i+1}-1\right)\\
h_i&\longrightarrow h_i-1
\ea
\right\}\text{with rate}\ q(\om_i,\,\om_{i+1}).\label{eq:rem}
\ee
Conditionally on the present state, these moves happen independently at all sites $i$. 
We can summarize this information in the  formal infinitesimal generator $L$
of the process  $\un\om(\cdot)$:
\be
\ba
(L\vp)(\un\om)&=\sum_{i\in\Zb}p(\om_i,\,\om_{i+1})\cdot\left[\vp(\dots,\,\om_i-1,\,\om_{i+1}+1,\,\dots)-\vp(\un\om)\right]\\
&+\sum_{i\in\Zb}q(\om_i,\,\om_{i+1})\cdot\left[\vp(\dots,\,\om_i+1,\,\om_{i+1}-1,\,\dots)-\vp(\un\om)\right].
\ea\label{eq:gen}
\ee
$L$ acts on bounded cylinder functions  $\vp\,:\,\Omega\to\Rb$ 
(this means that $\vp$ depends only on finitely many  $\om_i$-values). 

Thus we have a Markov process $\{\un\om(t): t\in\Rb^+\}$ of an evolving increment configuration and a 
Markov process
$\{\un h(t): t\in\Rb^+\}$ of an evolving height configuration. 
 The initial increments $\un\om(0)$ 
specify the initial height $\un h(0)$ up to a vertical translation.  We shall always
normalize the height process so that   \(h_0(0)=0\).

\begin{figure}[ht]
\begin{center}
\begin{pspicture}(160,160)
\psline{->}(0,5)(160,5)
\multirput(10,3)(20,0){8}{\psline(0,4)}
\multido{\n=-3+1}{8}{%
\setcounter{z}{70+20*\n}
\rput(\thez,0){\tiny\n}
}
\uput[315](160,5){\tiny\(i\)}
\rput(10,15){\(\bullet\)}
\rput(10,20){\(\bullet\)}
\rput(50,15){\(\bullet\)}
\rput(50,20){\(\bullet\)}
\rput(50,25){\(\bullet\)}
\rput(70,15){\(\bullet\)}
\rput(150,15){\(\bullet\)}
\psline[linewidth=1pt](0,160)(10,160)(10,120)(50,120)(50,60)(70,60)(70,40)(150,40)(150,30)
\psline[linewidth=1pt,linestyle=dashed,dash=1pt 2pt](50,80)(70,80)(70,60)
\psline[linestyle=dashed,dash=1pt 2pt]{->}(54,24)(66,21)
\psline[linestyle=dashed,dash=1pt 2pt]{->}(60,63)(60,77)

\psline[linestyle=dashed,dash=0.5pt 2pt](0,140)(10,140)
\psline[linestyle=dashed,dash=0.5pt 2pt](0,120)(10,120)
\psline[linestyle=dashed,dash=0.5pt 2pt](0,100)(50,100)
\psline[linestyle=dashed,dash=0.5pt 2pt](0,80)(50,80)
\psline[linestyle=dashed,dash=0.5pt 2pt](0,60)(50,60)
\psline[linestyle=dashed,dash=0.5pt 2pt](0,40)(70,40)

\psline[linestyle=dashed,dash=0.5pt 2pt](10,30)(10,120)
\psline[linestyle=dashed,dash=0.5pt 2pt](30,30)(30,120)
\psline[linestyle=dashed,dash=0.5pt 2pt](50,30)(50,60)
\psline[linestyle=dashed,dash=0.5pt 2pt](70,30)(70,40)
\psline[linestyle=dashed,dash=0.5pt 2pt](90,30)(90,40)
\psline[linestyle=dashed,dash=0.5pt 2pt](110,30)(110,40)
\psline[linestyle=dashed,dash=0.5pt 2pt](130,30)(130,40)
\psline[linestyle=dashed,dash=0.5pt 2pt](150,30)(150,40)

\rput[l](11,140){\(\Biggr\}\scriptstyle{\om_{-3}=2}\)}

\end{pspicture}
\caption{The wall and the particles with a possible step}\label{fig:wall}
\end{center}
\end{figure}
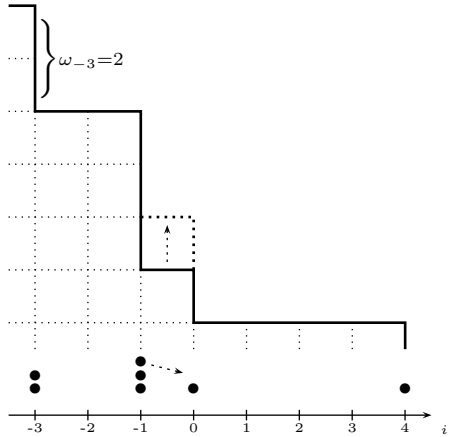

In the particle picture the variable \(\om_i(t)\) represents the number of particles 
at site $i$ 
at time \(t\).  Step \eqref{eq:add} represents a rightward jump of a particle over the edge \((i,\,i+1)\), while step \eqref{eq:rem} represents a leftward jump. (If negative \(\om\)-values 
are permitted, one needs to consider particles and  antiparticles, with antiparticles 
jumping in the opposite direction.) Figure \ref{fig:wall} shows a configuration and a possible step with both walls and particles. It is in the particle guise that many of these processes appear in the literature: simple exclusion processes, zero range processes and misanthrope processes are examples included in the class studied in this paper.

It will be useful to see that
\be
\ba
h_i(t)=h_i(t)-&h_0(0)=\text{the net number of particles that have passed, }\\
&\qquad\text{from left to right,  the straight-line space-time path}\\
&\qquad\text{that connects }(1/2,\,0)\text{ to }(i+1/2,\,t).
\ea\label{eq:hno}
\ee
In particular, height increment $h_i(t)-h_i(0)$ is the cumulative net particle current across the edge $(i,\,i+1)$ during time $(0,\,t]$. 

We impose the following four assumptions \eqref{eq:ominmax}--\eqref{eq:symm} on the rates.
\begin{itemize}
\item The rates $p,\,q:I\times I\to\Rb^+$  must satisfy
\be
p(\omin,\,\cdot\,)\equiv p(\,\cdot\,,\,\omax)\equiv q(\omax,\,\cdot\,)\equiv q(\,\cdot\,,\,\omin)\equiv0\label{eq:ominmax}
\ee
whenever either $\omin$ or $\omax$ is finite. Either both $p$ and $q$ are 
strictly positive in all other cases, or one of them is identically zero. 
The process is called  \emph{totally asymmetric} if either $q\equiv 0$ or $p\equiv 0$.
\item The dynamics has a smoothing effect when we assume 
 the following monotonicity:
\be
\ba
p(z+1,\,y)&\ge p(z,\,y),\qquad&p(y,\,z+1)&\le p(y,\,z)\\
q(z+1,\,y)&\le q(z,\,y),\qquad&q(y,\,z+1)&\ge q(y,\,z)
\ea\label{eq:mon}
\ee
for $y,\,z,\,z+1\in I$. Under this property  the higher the neighbors of a column, 
the faster it grows and the longer it waits for a brick removal, on average. 
This is the notion of  {\sl attractivity}.
\item The next two assumptions guarantee the existence of translation-invariant product-form
  stationary measures. (Similar assumptions were employed by  Cocozza-Thi\-vent \cite{coco}.)
\begin{itemize}
\item For any $x,\,y,\,z\in I$
\be
\ba
&p(x,\,y)+p(y,\,z)+p(z,\,x)\!\!\!\!\!&&\\
+\,&q(x,\,y)+q(y,\,z)+q(z,\,x)\!\!\!\!\!&=&\,p(x,\,z)+p(z,\,y)+p(y,\,x)\\
&&+&\,q(x,\,z)+q(z,\,y)+q(y,\,x).
\ea\label{eq:stacifelt}
\ee
\item There are symmetric functions $s_p$ and $s_q$ on $I\times I$, and a function $f$ 
on $I$ such that $f(\omin)=0$ whenever $\omin$ is finite, $f(z)>0$ for $z>\omin$,  
and for any $y,\,z\in I$, 
\be\begin{aligned}
p(y,\,z)&=s_p(y,\,z+1) f(y)\\
\text{and}\quad q(y,\,z)&=s_q(y+1,\,z) f(z).
\end{aligned}\label{eq:symm}\ee
(Interpret $s_p(y,\,z)=s_q(y,\,z)=0$ if $y$ or $z>\omax$.) 
Condition \eqref{eq:mon} implies that $f$ is nondecreasing on $I$.
\end{itemize}
\end{itemize}

An attempt at covering this broad class of processes  raises the uncomfortable point 
that there is no unified existence proof for this entire class. Different  constructions
 in the literature  place various boundedness or growth conditions
on $p$ and $q$ and the space $I$, and result in various degrees of regularity
 for the semigroup.  (Among key references are Liggett's monograph \cite{ips}, 
 and articles  \cite{and},  \cite{exists} and  \cite{lize}.)  
  These existence matters are beyond the scope of
  this paper.  Yet we 
wish to give a general proof for fluctuations that in principal works for all processes in the family,
subject to the more serious assumptions we explain in Section \ref{sc:microconc}. 
To avoid extraneous technical issues  we make the following blanket assumptions on the
rates $p$ and $q$ to be considered. 
 \begin{itemize}
 \item We assume that the increment process $\un\om(t)$, and the corresponding height 
 process $\un h(t)$ with normalization $h_0(0)=0$,  that obey 
 Poisson  rates $p$ and $q$ as described by \eqref{eq:add} and \eqref{eq:rem}, 
  can be constructed  with cadlag paths in a subspace 
 $\wt\Omega$ of tempered increment configurations (i.e.\ configurations that obey some restrictive growth conditions). 
\item The subspace 
  $\wt\Omega$ is of full measure under the invariant distributions $\mu^\te$ defined in Section \ref{sc:gibbs}. 
  \item  It is also possible to construct jointly several versions of the process
  with initial configurations from the space $\wt\Omega$ and with joint evolution obeying
  basic coupling (described in Section \ref{sc:basiccoupling}). 
\item  Rates  $p$ and $q$ have all moments under the invariant distributions $\mu^\te$. In fact arguments like Lemma \ref{lm:hreg} of the Appendix provide this when \(f\) does not grow faster than exponential on \(\Zb^+\) and does not decrease faster to zero than exponential on \(\Zb^-\).
\end{itemize}
 
 The reader will see that our proofs in Sections \ref{sc:ub}, \ref{sc:lb}, \ref{sc:Qlln} and 
 \ref{sc:zrp} do not make any analytic demands on the semigroup and its
 relation to the generator.  We only use couplings, 
counting of particle currents and simple Poisson bounds.

Two identities from article \cite{varj2nd} play a key role in this paper, given
as \eqref{eq:varcovar} and \eqref{eq:qchar} in Section \ref{sc:hydro}. 
 These identities
hold for all processes in the family under study.   The proofs   given in \cite{varj2nd} use
generator calculations which may not be justified for all these processes.
However, these identities can also be proved by counting particles and 
taking limits of finite-volume processes (\cite{asepsimple} contains an example).
Such a proof should be available with any reasonable construction of a process. 
Hence we shall not hesitate to use the results of \cite{varj2nd}.

\subsection{Examples}\label{sc:genex}

To give concrete meaning to the general formulation  of the previous section we describe
some basic examples. 
The type of state space $I$ distinguishes three cases that we call generalized exclusion, 
misanthrope and bricklayers processes.  In all cases there are two parameters
  $0\le p,\,q\le 1$ such that  $p+q=1$. \emph{Asymmetric processes} have $p\ne q$.
  These are the processes for which our results are relevant. 

\begin{enumerate}
\item {\bf Generalized exclusion processes.} These are  the cases where 
 both $\omin$ and $\omax$ are finite.
\begin{itemize}
\item {\bf The asymmetric simple exclusion process (ASEP)} introduced by F.\ Spitzer \cite{spi} is defined by $\omin=0,\ \omax=1$, $f(z)={\bf1}\{z=1\}$,
$
s_p(y,\,z)=p\cdot{\bf1}\{y=z=1\}$ and $s_q(y,\,z)=q\cdot{\bf1}\{y=z=1\}$. 
This produces the familiar rates 
\[
p(y,\,z)=p\cdot{\bf1}\{y=1,\,z=0\}\quad\text{and}\quad q(y,\,z)=q\cdot{\bf1}\{y=0,\,z=1\}.
\]
Here $\om_i\in\{0,\,1\}$ is the occupation number for site $i$, $p(\om_i,\,\om_{i+1})$ is the rate for a particle to jump from site $i$ to $i+1$, and $q(\om_i,\,\om_{i+1})$ is the rate for a particle to jump from site $i+1$ to $i$. These rates have values $p$ and $q$, respectively, whenever there is a particle to perform the above jumps, and there is no particle on the terminal site of the jumps. Conditions \eqref{eq:mon} and \eqref{eq:stacifelt} are also satisfied by these rates.
\item {\bf Particle-antiparticle exclusion process.}   Let   $\omin=-1$, $\omax=1$. 
Take $f(-1)=0$, $f(0)=c$ ({\sl creation}), $f(1)=a$ ({\sl annihilation})
where $c$ and $a$ are positive rates with $c\le a/2$,
\[
\ba
s_p(0,\,1)=s_p(1,\,0)&=p,&\ s_p(0,\,0)&=\frac{pa}{2c},&\ s_p(1,\,1)&=\frac{p}{2},\\
s_q(0,\,1)=s_q(1,\,0)&=q,&\ s_q(0,\,0)&=\frac{qa}{2c},&\ s_q(1,\,1)&=\frac{q}{2}
\ea
\]
and $s_p,\ s_q$ zero in all other cases. These result in rates
\[
\ba
p(0,\,0)&=pc,&\ p(0,\,-1)=p(1,\,0)&=\frac{pa}{2},&\ p(1,\,-1)&=pa,\\
q(0,\,0)&=qc,&\ q(-1,\,0)=q(0,\,1)&=\frac{qa}{2},&\ q(-1,\,1)&=qa
\ea
\]
and zero in all other cases. If $\om_i$ is the number of particles at site $i$, with $\om_i=-1$ meaning the presence of an antiparticle, then this model describes an asymmetric exclusion process of particles and antiparticles with annihilation and particle-antiparticle pair creation. These rates also satisfy our conditions.
\end{itemize}
One can imagine 
other generalizations with  bounded numbers of particles and/or antiparticles per site.
\item {\bf Generalized misanthrope processes} have $\omin>-\infty,\ \omax=\infty$.
\begin{itemize}
\item {\bf Zero range process.} Take  $\omin=0,\ \omax=\infty$, an arbitrary nondecreasing function $f\,:\,\Zb^+\to\Rb^+$ such that $f(0)=0$,
\[
\ba
s_p(y,\,z)&\equiv p\qquad&&\text{and}&\qquad s_q(y,\,z)&\equiv q,\\
p(y,\,z)&=p f(y)\qquad&&\text{and}&\qquad q(y,\,z)&=q f(z).
\ea
\]
Again, $\om_i$ represents the number of particles at site $i$. Depending on this number, a particle jumps from $i$ to the right with rate $p f(\om_i)$, and to the left with rate $q f(\om_i)$. 
These rates trivially satisfy conditions \eqref{eq:mon} and \eqref{eq:stacifelt}.
\end{itemize}
\item {\bf General deposition processes} have $\omin=-\infty$ and $\omax=\infty$. The height differences between adjacent columns  can be arbitrary integers. Antiparticles 
are needed for a particle representation of the process.
\begin{itemize}
\item {\bf Bricklayers process.} Let $f\,:\,\Zb\to\Rb^+$ be non-decreasing and satisfy \[
f(z)\cdot f(1-z)=1\qquad\text{for all }z\in\Zb.
\]
The values of $f$ for positive $z$'s thus determine the values for non-positive $z$'s. Let 
\[
s_p(y,\,z)=p+\frac{p}{f(y)f(z)}\qquad\text{and}\qquad s_q(y,\,z)=q+\frac{q}{f(y)f(z)},
\]
which results in
\[
p(y,\,z)=pf(y)+pf(-z)\qquad\text{and}\qquad q(y,\,z)=qf(-y)+qf(z).
\]
The following picture motivates the name bricklayers process. At each site $i$ stands a 
  bricklayer who lays a brick on the column to his left at rate $pf(-\om_i)$ and on
  the column to his right at rate $pf(\om_i)$. 
  Each bricklayer  also removes a brick from his left at rate $qf(\om_i)$ and from 
  his right at rate $qf(-\om_i)$.   Conditions \eqref{eq:mon} and \eqref{eq:stacifelt} hold for the rates.
\end{itemize}
\end{enumerate}

These were  examples for which our theorem holds, \emph{provided 
the hypotheses on microscopic concavity to be described below can be verified.}

\subsection{Basic coupling}
\label{sc:basiccoupling}
In  \emph{basic coupling} the joint evolution of \(n\) processes \(\un\om^m(\cdot),\ m=1,\dotsc, n\), is defined in such a manner that the processes ``jump together as much as possible.''
The joint rates are determined as follows, given the current configurations
  \(\un\om^1,\,\un\om^2,\,\dots,\,\un\om^n\in\wt\Omega\).  Consider  a step
    of type \eqref{eq:add} over the edge \((i,\,i+1)\).   Let $m\mapsto\ell(m)$ be a permutation that orders  the rates
  of the individual processes for this move:
\[
r(m)\equiv p(\om^{\ell(m)}_i, \, \om^{\ell(m)}_{i+1})\le p(\om^{\ell(m+1)}_i, \, \om^{\ell(m+1)}_{i+1})
\equiv r(m+1),\quad 1\le m< n.
\]
Set also the dummy value \( r(0)=0\).  
Now the rule is that independently 
  for each \(m=1,\dotsc, n\), 
  at rate $r(m)-r(m-1)$,  precisely processes 
   $\un\om^{\ell(m)}$, $\un\om^{\ell(m+1)}$, $\dots$, $\un\om^{\ell(n)}$
    execute move \eqref{eq:add}, and processes
  \(\un\om^{\ell(1)},\,\un\om^{\ell(2)},\,\dots,\,\un\om^{\ell(m-1)}\) do not. 
  The combined effect of these joint rates 
   creates the correct marginal rates, that is, process $\un\om^{\ell(m)}$ executes
  this move with rate $r(m)$. 
  
    Notice also that, due to \eqref{eq:mon}, a jump of \(\un\om^a\) without \(\un\om^b\) can only occur 
    if \(p(\om^b_i,\,\om^b_{i+1})<p(\om^a_i,\,\om^a_{i+1})\) which implies \(\om^a_i>\om^b_i\) or \(\om^a_{i+1}<\om^b_{i+1}\). The result of this step \eqref{eq:add} then cannot increase the 
    number of discrepancies  between the two processes, hence the name \emph{attractivity} for \eqref{eq:mon}.  In particular, a sitewise ordering $\om^a_i\le \om^b_i$\quad\(\forall i\in\Zb\) is preserved by
    the basic coupling. 

One can check that  moves   of type \eqref{eq:rem} with rates  \(q\) obey the same attractivity property.

The differences between two processes are called \emph{second class particles}. Their number is nonincreasing. In particular, if \(\om^a_i\ge\om^b_i\) for each \(i\in\Zb\), then the second class particles are conserved. In view of \eqref{eq:hno}, in this case the net number of second class particles 
that pass from left to right across
 the straight-line space-time path from \((1/2,\,0)\) to \((i+1/2,\,t)\)  equals the growth difference
\be
\bigl(h^a_i(t)-h^a_0(0)\bigr)-\bigl(h^b_i(t)-h^b_0(0)\bigr)=h^a_i(t)-h^b_i(t)\label{eq:2hno}
\ee
between the two processes \(\un\om^a(\cdot)\) and \(\un\om^b(\cdot)\).

 A special case that is 
of key importance to us is the situation 
where only one second class particle is present between two processes.

\subsection{Translation invariant stationary product distributions}\label{sc:gibbs}

The results of this paper concern stationary processes with particular product-form
marginal distributions that we define in this section. 
  For many cases it has been proved that these measures  are the only extremal translation-invariant stationary distributions. Following some ideas in Cocozza-Thivent \cite{coco}, we first consider the nondecreasing function $f$ whose existence was assumed in  \eqref{eq:symm}.
 For $I\ni z>0$ define 
\[
f(z)!:\,=\prod_{y=1}^zf(y),
\]
while for $I\ni z<0$ let
\[
f(z)!:\,=\frac{1}{\prod\limits_{y=z+1}^0f(y)},
\]
and then  $f(0)!:\,=1$. This definition satisfies 
$
f(z)!\cdot f(z+1)=f(z+1)!
$
for all $z\in I$. Let
\[
\bar\te:\,=\left\{\begin{array}{ll}\log\left(\liminf\limits_{z\to\infty}\left(f(z)!\right)^{1/z}\right)=\lim\limits_{z\to\infty}\log(f(z)),\ \ &\ \text{if}\ \omax=\infty\\\infty,\ \ &\ \text{else}\end{array}\right.
\]
and
\[
\un\te:\,=\left\{\begin{array}{ll}\log\left(\limsup\limits_{z\to\infty}\left(f(-z)!\right)^{-1/z}\right)=\lim\limits_{z\to\infty}\log(f(-z)),\ \ &\ \text{if}\ \omin=-\infty\\-\infty,\ \ &\ \text{else}.\end{array}\right.
\]
By monotonicity of $f$, we have $\bar\te\ge\un\te$. The case \(\bar\te=\un\te\) would imply that \(\omin=-\infty,\ \omax=\infty\), and \(f\) is a constant. Notice that \eqref{eq:mon} and \eqref{eq:symm} imply that \(s_p\) is non-increasing in its variables, but \(p\) is non-decreasing in its first variable. Hence a constant \(f\) results in an \(s_p\) that does not depend on its first variable. But then by its symmetric property it does not depend on its second variable either, and we conclude that a constant \(f\) implies constant rates \(p\) (and, similarly, \(q\)). We exclude this uninteresting case by postulating
\be
\text{Assume \(f\) to be such that \(\un\te<\bar\te\).}
\label{eq:f-ass}\ee

 For $\te\in\left(\un\te,\,\bar\te\right)$ define the state sum
\be
Z(\te):\,=\sum_{z\in I}\frac{\e{\te z}}{f(z)!}<\infty.\label{eq:zdef}
\ee
Let the product-distribution $\un\mu^\te$ on $\Omega=I^{\Zb}$ have marginals
\be
\mu^\te(z)=\un\mu^\te\left\{\un\om\,:\,\om_i=z\right\}:\,=
\frac{1}{Z(\te)}\cdot\frac{\e{\te z}}{f(z)!}
\qquad(z\in I).
\label{eq:mudef}
\ee
Assumptions \eqref{eq:ominmax}, \eqref{eq:mon}, \eqref{eq:stacifelt}, \eqref{eq:symm} imply that
for $\te\in\left(\un\te,\,\bar\te\right)$ 
 the product distribution $\un\mu^\te$ is stationary for the process generated by \eqref{eq:gen} (see \cite{varj2nd}; notice that the top display on Page 443 of \cite{varj2nd} is incorrect, to get the correct identity, multiply with the cylinder functions and take expectation). For some calculations in the Appendix it will be convenient to 
note that the family  \(\{\mu^\te\}\) can be obtained by exponentially weighting
a probability measure  \(\mu^{\te_0}\) for a fixed value \(\te_0\in(\un\te,\,\bar\te)\).

\(\Pv^\te,\ \Ev^\te,\ \Vv^\te,\ \Cov^\te\) will refer to laws of a process evolving in this stationary distribution. In the appendix we show that the \emph{density}
\[
\vr(\te):\,=\Ev^\te(\om)
\]
is a strictly increasing, infinitely differentiable  function of the parameter \(\te\)
that maps the interval $(\un\te,\,\bar\te)$ onto the interval
$(\omin,\,\omax)$. (The following point should cause no confusion:  the single-site
state space $I$ consists of the integers between $\omin$ and $\omax$, including 
endpoints if finite,  but for density values the interval $(\omin,\,\omax)$ is
an interval of real numbers.) 
For most cases we shall use the density $\vr$, rather than \(\te\), for parameterizing the stationary distributions. Accordingly, \(\mu^\vr,\ \Pv^\vr,\ \Ev^\vr,\ \Vv^\vr,\ \Cov^\vr\) will refer to laws of a density \(\vr\) stationary process.

\subsection{Hydrodynamics and some exact identities}
\label{sc:hydro}
The \emph{hydrodynamic flux} is defined as
\be
\Hc(\vr):\,=\Ev^\vr(p(\om_0,\,\om_1)-q(\om_0,\,\om_1)).\label{eq:fluxdef}
\ee
$\Hc(\vr)$ is the expected net rate at which a given column grows, or at which particles pass any fixed 
lattice edge from left to right in a stationary density-$\vr$ process.  
We show smoothness of $\Hc$ in Section \ref{sc:afx} of the Appendix. It is expected, 
and in many instances proved, that asymmetric members of our class satisfy the conservation law
\[
\partial_T\vr(T,\,X)+\partial_X\Hc(\vr(T,\,X))=0
\]
in the Eulerian-scaled time and space variables \(T\) and \(X\), see e.g.\ Rezakhanlou \cite{hl} or Bahadoran, Guiol, Ravishankar and Saada \cite{bagurasa}. The \emph{characteristic speed} is the velocity with which small perturbations of the solution of this PDE propagate, and is given by
\be
V^\vr:\,=\Hc'(\vr).\label{eq:chardef}
\ee

A particular expectation we shall need several times
is 
\be
\Ev^\vr(h_i(t))=\Hc(\vr)t-\vr i, \quad t\ge 0,\; i\in\Zb.
\label{eq:Eh}\ee
For $i=0$ this follows from \eqref{eq:hno}, and in general from the $i=0$ case 
together with $\om_j(t)=h_{j-1}(t)-h_j(t)$.  

When a stationary process is perturbed  by adding a second class particle at the
origin at time zero, we obtain two processes, \(\uomm(\cdot)\) and  \(\uomp(\cdot)\). It is
not a priori clear what the initial joint distribution of the occupation variables $\omm_0(0)$, \(\om_0(0)\)
should be.  For ASEP there is no ambiguity due to the simplicity
of the single-site state space: the only way to have  a discrepancy is to set
 $\omm_0(0)=0$, $\omp_0(0)=1$.  A useful generalization  of this
distribution to the broader class of processes 
 involves  the following  family of probability measures on $I$  
 introduced in \cite{varj2nd}:
\be
\wih\mu^\vr(y):\,=\frac{1}{\Vv^\vr(\om_0)}\sum_{z=y+1}^{\omax}(z-\vr)\mu^\vr(z),\qquad y\in I.\label{eq:muhat}
\ee
 An empty sum is zero by convention and so if $\omax<\infty$, $\wih\mu^\vr(\omax)=0$.
Consequently there is room for an additional particle under the $\wih\mu^\vr$ distribution,
in the sense that if $\om\sim\wih\mu^\vr$ then also $\om+1\in I$.

 To our knowledge these distributions $\wih\mu^\vr$
 do not possess any invariance properties. 
 Their virtue is that  they make  identities \eqref{eq:varcovar} and \eqref{eq:qchar} below true. 
 We show in Section \ref{sc:cvxi} of the Appendix that both \(\mu^\vr\) and \(\wih\mu^\vr\) are stochastically monotone in the density \(\vr\). 
 (There is, however, no stochastic domination
   between \(\mu^\vr\) and \(\wih\mu^\vr\) in general.)

Denote by \(\Ev\) the expectation w.r.t.\ the evolution of a pair \((\uomm(\cdot),\,\uomp(\cdot))\) started with initial data (recall \eqref{eq:dedef})
\be
\uomm(0)=\uomp(0)-\un\de_0\sim\Bigl(\bigotimes_{i\ne0}\mu^\vr\Bigr)\otimes\wih\mu^\vr,\label{eq:umhdef}
\ee
and evolving under the basic coupling.  
 This pair will always have a single second class 
particle whose position is denoted by \(Q(t)\).  In other words, 
$\uomm(t)=\uomp(t)-\un\de_{Q(t)}$. 
Corollaries 2.4 and 2.5 of \cite{varj2nd} state that 
\begin{align}
\Vv^\vr(h_i(t))&=\Vv^\vr(\om)\cdot\Ev|Q(t)-i|\label{eq:varcovar}
\intertext{and}
\Ev(Q(t))&=V^\vr\cdot t\label{eq:qchar}
\end{align}
for any \(i\in\Zb\) and \(t\ge0\).  Note in particular that in \eqref{eq:varcovar} the variances
are taken in a stationary process, while the expectation of $Q(t)$ is taken in the
coupling with initial distribution \eqref{eq:umhdef}. 
These two identities follow from the definition of our models together with translation invariance and the product structure of the stationary distribution.
 
\subsection{Microscopic concavity}
\label{sc:microconc} 
From now on fix the jump rates $p, q:I\times I\to\Rb^+$ that define the process
in question, assumed to satisfy all the assumptions discussed thus far. 
The \(t^{1/3}\) current or height fluctuations are expected 
when the hydrodynamic flux \(\Hc(\vr)\) is strictly concave or convex. In this paper we  
 discuss only  the strictly concave case.  This implies that the characteristic speed $V^\vr=\Hc'(\vr)$ 
is a decreasing function of density $\vr$:
\be
\la<\vr\;\Longrightarrow V^\la>V^\vr.  \label{eq:Vorder}\ee

 The microscopic counterpart 
of  a characteristic is the motion of a second class particle.  Our key assumption
that we term \emph{microscopic concavity} is that the ordering \eqref{eq:Vorder} 
can also be realized at the particle level as an ordering between two second class 
particles introduced into two processes at densities $\la$ and $\vr$.  Since this 
is now a probabilistic notion, there are several possible formulations, ranging from
almost sure ($Q^\la(t)\ge Q^\vr(t)$ in a coupling) to distributional formulations. 
Assumption \ref{def:microc} below gives the precise technical form in which
this paper utilizes this notion of microscopic concavity.  It stipulates that the
ordering of second class particles is achieved by processes that 
evolve on the labels of auxiliary second class particles, and also requires some control
of the tails of these random labels.  

We do not imagine that this precise 
formulation will be the right one for all processes. We take it as a starting point
and  future work may lead to alternative formulations.
Assumption \ref{def:microc} has the virtue that its requirements  can be verified for some interesting processes.

 Let \(\la<\vr\) be two densities.
 Proposition \ref{pr:muhatpr} in the Appendix
gives the stochastic domination  \({\wih\mu}^\la\le{\wih\mu}^\vr\).
Define \({\wih\mu}^\vr+1\) as the measure that gives weight
\({\wih\mu}^\vr(z-1)\) to an integer $z$ such that  \(\omin<z<\omax+1\).
Let \({\wih\mu}^{\la,\vr}\) be a coupling measure with
marginals \({\wih\mu}^\la\) and \({\wih\mu}^\vr+1\) and with the property
\be
{\wih\mu}^{\la,\vr}\{(y,\,z)\,:\,\omin-1<y<z<\omax+1\}=1.\label{eq:cmuhat}
\ee
Let also
\(\mu^{\la,\vr}\) be a coupling measure of site-marginals \(\mu^\la\)
and \(\mu^\vr\) of the invariant distributions,  with
\be
\mu^{\la,\vr}\{(y,\,z)\,:\,\omin-1<y\le z<\omax+1\}=1,\label{eq:cmu}
\ee
this is possible by Corollary \ref{cr:must} of the Appendix. Note the distinction that under ${\wih\mu}^{\la,\vr}$ the second coordinate
is strictly above the first. 

To have notation for inhomogeneous product measures on $I^\Zb$, let 
\(\un\la=(\la_i)_{i\in\Zb}\) and \(\un\vr=(\vr_i)_{i\in\Zb}\) denote sequences
of density values, with  \(\la_i\) and \(\vr_i\) assigned to site \(i\).  
The product distribution with marginals \({\wih\mu}^{\la_0,\vr_0}\)
at the origin and \(\mu^{\la_i,\vr_i}\) at other sites is denoted by
\be
{\un{\wih\mu}}^{\un\la,\un\vr}:\,=\Bigl(\bigotimes_{i\ne0}\mu^{\la_i,\vr_i}\Bigr)\otimes{\wih\mu}^{\la_0,\vr_0}.\label{eq:unmu}
\ee
Measure \({\un{\wih\mu}}^{\un\la,\un\vr}\) gives probability one to the event
\[
\{(\un\eta(0),\,\uomp(0))\,:\,\eta_0(0)<\omp_0(0)\text{, and }\eta_i(0)
\le\omp_i(0)\text{ for }0\ne i\in\Zb\}.
\]
The initial configuration \((\un\eta(0),\,\uomp(0))\) will always be assumed a member
of this set, and the pair process \((\un\eta(t),\,\uomp(t))\) evolves in basic coupling. 
In general  \({\un{\wih\mu}}^{\un\la,\un\vr}\) is \emph{not} stationary for this joint evolution.

The discrepancies between these two processes are called the
\emph{ \(\omp-\eta\) (second class) particles}. The number of such particles at site
 $i$ at time $t$ is $\omp_i(t)-\eta_i(t)$.   In the basic coupling   the
 \(\omp-\eta\) particles are
conserved, in the sense that none are created or annihilated.
 We label the  \(\omp-\eta\) particles with
 integers, and let $X_m(t)$ denote the position of particle $m$ at time $t$. 
 The  initial 
labeling is chosen to   satisfy
\[
\dotsm\le X_{-1}(0)\le X_0(0)=0<X_{1}(0)\le\dotsm.
\]
We can specify that $X_0(0)=0$ because under  \({\un{\wih\mu}}^{\un\la,\un\vr}\)
there is an \(\omp-\eta\) particle at site $0$ with probability $1$. 
During the evolution we keep the positions $X_i(t)$  of the  \(\omp-\eta\) particles ordered.
To achieve this we stipulate that 
\be\begin{array}{l}
\text{whenever an  \(\omp-\eta\) particle jumps from
a site, }\\
\text{if the jump is to the right  the highest label moves,}\\
\text{and if the jump is to the left  the  lowest label moves.}\end{array}\label{eq:Xrule}\ee

Here is the precise form of microscopic concavity for this paper.  The assumption 
states that a certain joint construction of processes (that is, a coupling) can be
performed for a range of densities in a neighborhood of a fixed density $\vr$. 
Recall \eqref{eq:dedef} for the definition of the configuration \(\un\de\).

\begin{ass}\label{def:microc}
Given a density $\vr\in(\omin,\,\omax)$, there exists $\ga_0>0$ such that the  following holds. 
For any  $\un\la$ and $\un\vr$ such that  \(\vr-\ga_0\le \la_i\le\vr_i\le\vr+\ga_0\) for all \(i\in\Zb\), 
 a joint process  \((\un\eta(t),\,\uomp(t),\,y(t),\,z(t))_{t\ge 0}\) can be constructed
with the following properties.  
\begin{itemize} 
\item Initially \((\un\eta(0),\,\uomp(0))\) is \({\un{\wih\mu}}^{\un\la,\un\vr}\)-distributed
and the joint process  \((\un\eta(\cdot),\,\uomp(\cdot))\) evolves in basic coupling.
\item Processes $y(\cdot)$ and  $z(\cdot)$ are integer-valued. Initially \(y(0)=z(0)=0\).
  With probability one
\be
y(t)\le z(t)\text{ for all }t\ge0.\label{eq:yz}
\ee
\item Define the processes 
\be
\uomm(t):\,=\uomp(t)-\un\de_{X_{y(t)}(t)}\quad\text{and}\quad\uetap(t):\,=\un\eta(t)+\un\de_{X_{z(t)}(t)}.\label{eq:indef}
\ee
Then both pairs $(\un\eta,\,\uetap)$ and  $(\uomm,\,\uomp)$  evolve marginally in basic
coupling. 
\item  For each $\ga\in(0,\,\ga_0)$ and large enough $t\ge 0$ 
 there exists a probability  distribution \(\nu^{\vr,\ga}(t)\) 
on \(\Zb^+\)  satisfying the tail bound
\be
\nu^{\vr,\ga}(t)\{y\,:\,y\ge y_0\}\le Ct^{\kappa-1}\ga^{2\kappa-3}y_0^{-\kappa}\label{eq:nuc}
\ee
for some fixed constants  \(3/2\le\kappa<3\) and $C<\infty$, and 
 such that  if   \(\vr-\ga\le \la_i\le\vr_i\le\vr+\ga\) for all \(i\in\Zb\), 
 then we have the stochastic bounds 
 \be y(t)\overset{\text{d}}{\le}\nu^{\vr,\ga}(t)\quad\text{and}\quad 
z(t)\overset{\text{d}}{\ge}-\nu^{\vr,\ga}(t).  \label{eq:yzbds}\ee
\end{itemize}
\end{ass}

Let us clarify some of the details in this assumption. 

Equation \eqref{eq:indef} says that \(Q^\eta(t):\,=X_{z(t)}(t)\) is  the single
second class particle between \(\un\eta\) and \(\uetap\),
 while \(Q(t):\,=X_{y(t)}(t)\) is  the one between \(\uomm\) and \(\uomp\). 
The first three bullets   say that it is possible to construct jointly four processes 
 \((\un\eta,\,\uetap,\,\uomm,\,\uomp)\) with the  specified initial conditions and so that 
 each pair 
  \((\un\eta,\,\uomp)\), $(\un\eta,\,\uetap)$ and $(\uomm,\,\uomp)$ has
   the desired marginal distribution,
  and most importantly so that 
\be
Q^\eta(t)=X_{z(t)}(t)\ge X_{y(t)}(t)=Q(t).\label{eq:qq}
\ee
 This is a consequence of \eqref{eq:yz}  because the 
  \(\omp-\eta\) particles $X_i(t)$  stay ordered.

The tail bound \eqref{eq:nuc} is formulated in this somewhat complicated fashion because
this appears  to be the weakest form our present proof allows.  In our currently available examples 
$\nu^{\vr,\ga}(t)$ is actually a fixed geometric distribution.  However, we expect that 
other examples will require more complicated bounds and so including this generality
is sensible. 

The assumptions made imply \(\un\eta(t)\le\uomp(t)\) a.s.,
and by   \eqref{eq:indef}  
\[
\un\eta(t)\le\uetap(t)\le\uomp(t)\quad\text{and}\quad\un\eta(t)\le\uomm(t)\le\uomp(t)\quad\text{a.s.}
\]
In our actual constructions of the processes \(\un\eta,\,\uetap,\,\uomm,\,\uomp\) for ASEP (Section \ref{sc:asep} and \cite{asepsimple}), for a class of totally asymmetric zero range processes (Section \ref{sc:zrp}) and for the totally asymmetric bricklayers process with exponential rates (future work) it turns out that the triples \((\un\eta,\,\uetap,\,\uomp)\) and \((\un\eta,\,\uomm,\,\uomp)\) evolve also in
 basic coupling, but the full joint evolution \((\un\eta,\,\uetap,\,\uomm,\,\uomp)\) does not.

As already explained, the microscopic concavity idea is contained in 
inequality \eqref{eq:yz}. There is also a sense in which 
the tail bounds \eqref{eq:yzbds} relate to concavity of the flux. 
Consider the situation
  \(\la_i\equiv\la<\vr\equiv\vr_i\).  We would expect the \(\omp-\eta\) particle \(X_0(\cdot)\)
to have average and long-term velocity
\[
R(\la,\,\vr)=\frac{\Hc(\vr)-\Hc(\la)}{\vr-\la},
\]
the Rankine-Hugoniot or shock speed. By concavity 
 \(\Hc'(\vr)=V^\vr\le R(\la,\,\vr)\le V^\la=\Hc'(\la)\). A strict microscopic counterpart 
 would be \(y(t)\le0\le z(t)\). But this condition is overly restrictive.  The only cases we know
 to satisfy it are  the totally asymmetric simple exclusion process and the
 totally asymmetric zero range process with  constant rate.  The 
 distributional bounds \eqref{eq:yzbds} are natural relaxations of  \(y(t)\le0\le z(t)\).

 By the same token,  perhaps the way to covering more examples with our
 approach  involves
 a similar distributional weakening of  \eqref{eq:yz}, but this seems less straightforward.
  
\subsection{Results}

We need a few more assumptions
and  then we can state the main result. Constants \(C_{\centerdot},\,\al_{\centerdot}\) will not depend on time, but might depend on the density parameter \(\vr\), and their values can change from line to line.  We are now working with  a fixed member of the class of processes described in
Section \ref{sc:modfam} with rate functions $p,\,q:I\times I\to\Rb^+$.  
 Recall that $\Hc$ is the hydrodynamic flux  defined in \eqref{eq:fluxdef}. In the Appendix
we show  $\Hc$ is infinitely differentiable under the restrictions on the rates placed in Section
\ref{sc:modfam}.  

\begin{ass}\label{ass:main}  The rates $p,\,q$ and 
density \(\vr\in(\omin,\,\omax)\) have  the following properties. 
\begin{itemize}
\item  The jump rate functions $p$ and $q$  satisfy assumptions \eqref{eq:ominmax}, \eqref{eq:mon}, \eqref{eq:stacifelt}, \eqref{eq:symm} and \eqref{eq:f-ass} discussed in Sections \ref{sc:modfam} 
and \ref{sc:gibbs}.
\item   \(\Hc''(\vr)<0\). 
\item Let $(\uomm,\,\uomp)$ be a pair of processes in basic coupling,
started from distribution \eqref{eq:umhdef}, with second class particle $Q(t)$. 
Then there exist constants \(0<\al_0,\,C_0<\infty\) such that 
\be
\Pv\{|Q(t)|>K\}\le
C_0\cdot\frac{t^2}{K^3}
\label{eq:expq}
\ee
whenever \(K>\al_0 t\) and $t$ is large enough.
\end{itemize}
\end{ass}

As mentioned, our results are valid only for asymmetric processes. The assumption
of asymmetry is implicitly contained in $\Hc''(\vr)<0$.  Symmetric processes 
have $\Hc(\vr)\equiv 0$.  
Exponential tail bounds for $|Q(t)|$ that imply assumption \eqref{eq:expq} hold automatically if the rates $p,\,q$ 
have bounded increments
because the rates for $Q$ come from these increments of $p$ and $q$. 
Here is the main result. 

\begin{tm}\label{tm:main}
Let  Assumptions \ref{def:microc} and \ref{ass:main} hold for density $\vr$. 
Let the processes $(\uomm(t),\,\uomp(t))$ evolve in basic coupling with
initial distribution {\rm\eqref{eq:umhdef}} and let $Q(t)$ be the position of the second class
particle between $\uomm(t)$ and $\uomp(t)$.
Then there is a constant \(C_1=C_1(\vr)\in(0,\,\infty)\) such that for all \(1\le m<3\),
\be
\frac1{C_1}<\liminf_{t\to\infty}\frac{\Ev|Q(t)-V^\vr t|^m}{t^{2m/3}}\le\limsup_{t\to\infty}\frac{\Ev|Q(t)-V^\vr t|^m}{t^{2m/3}}<\frac{C_1}{3-m}.
\label{eq:main}\ee
\end{tm}
Diffusive fluctuations are characterized by a variance of order \(t\). The estimates above show that the second class particle has variance of order \(t^{4/3}\), this is called superdiffusivity.

Next some corollaries. 
Notation \(\fl{X}\) stands for the lower integer part of \(X\).
\begin{cor}[Current variance]\label{cr:varj}
Under Assumptions \ref{def:microc} and  \ref{ass:main}, there is a constant \(C_1=C_1(\vr)>0\), such that
\[
\frac1{C_1}<\liminf_{t\to\infty}\frac{\Vv^\vr(h_{\fl{V^\vr t}}(t))}{t^{2/3}}\le\limsup_{t\to\infty}\frac{\Vv^\vr(h_{\fl{V^\vr t}}(t))}{t^{2/3}}<C_1.
\]
\end{cor}
\noindent
This follows from \eqref{eq:varcovar} with the choice \(m=1\).
\begin{cor}[Law of Large Numbers for the second class particle]\label{cr:lln}
Under Assumptions \ref{def:microc} and  \ref{ass:main}, the Weak Law of Large Numbers holds 
in a density-\(\vr\) stationary process:
\be
\frac{Q(t)}{t}\overset{\text{d}}{\to}V^\vr.
\label{eq:Qlln}\ee
If the rates \(p\) and \(q\) 
have bounded increments, then  almost sure
convergence also holds in {\rm\eqref{eq:Qlln}} (Strong Law of Large Numbers).
\end{cor}
The Weak Law is a simple consequence of Theorem \ref{tm:main}. 
The Strong Law will be proved in Section \ref{sc:Qlln}.
\begin{cor}[Dependence of current on the initial configuration]
Under Assumptions \ref{def:microc} and  \ref{ass:main}, for any \(V\in\Rb\) and \(\al>1/3\) the following limit holds in the \(L^2\) sense for a density-\(\vr\) stationary process:
\be
\lim_{t\to\infty}\frac{h_{\fl{Vt}}(t)-h_{\fl{Vt}-\fl{V^\vr t}}(0)-t(\Hc(\vr)-\vr\Hc'(\vr))}{t^\al}=0.
\label{eq:L2cor}\ee
\label{cr:initial}\end{cor}
Recall that
\be
h_{\fl{Vt}-\fl{V^\vr t}}(0)=\left\{\ba
&\sum_{i=\fl{Vt}-\fl{V^\vr t}+1}^0\om_i(0),&&\text{if }V<V^\vr,\\
&\qquad\quad0,&&\text{if }V=V^\vr,\\
&-\sum_{i=1}^{\fl{Vt}-\fl{V^\vr t}}\om_i(0),&&\text{if }V>V^\vr
\ea\right.
\label{eq:inith}\ee
only depends on a finite segment of the initial configuration. Limit
\eqref{eq:L2cor} shows that on the 
diffusive time scale $t^{1/2}$  only fluctuations from the initial distribution are visible: these
fluctuations are  translated rigidly at the characteristic speed \(V^\vr\). 
Proof of \eqref{eq:L2cor} follows by translating 
\(h_{\fl{Vt}}(t)-h_{\fl{Vt}-\fl{V^\vr t}}(0)\) to \(h_{\fl{V^\vr t}}(t)-h_0(0)=h_{\fl{V^\vr t}}(t)\) and by applying Corollary \ref{cr:varj}. From \eqref{eq:L2cor}, \eqref{eq:inith} and the i.i.d.\ initial
$\{\om_i\}$ follow a limit for the variance and a
central limit theorem (CLT), which we record  in our final corollary. 
\begin{cor}[Central Limit Theorem for the current]
Under Assumptions \ref{def:microc} and  \ref{ass:main}, for any \(V\in\Rb\) in a density-\(\vr\) stationary process
\be
\lim_{t\to\infty}\frac{\Vv^\vr(h_{\fl{Vt}}(t))}{t}=\Vv^\vr(\om)\cdot|V^\vr-V|=\,:D,\label{eq:normal}
\ee
and the Central Limit Theorem also holds: the centered and normalized  height
 \(\wt h_{\fl{Vt}}(t)/\sqrt{t\cdot D}\) converges in distribution to a standard normal.
\label{cr:CLT}\end{cor}
\noindent

For ASEP the CLT, the limiting variance 
  \eqref{eq:normal} and the appearance  of initial fluctuations on the diffusive scale 
  were proved   by P.\ A.\ Ferrari and L.\ R.\ G.\ Fontes \cite{se}.
 For convex rate zero range and bricklayers processes Corollary \ref{cr:CLT}
 was proved  by M.\ Bal\'azs \cite{fluct}.  
   
\subsubsection*{Remark on the convex case}
 Our results and proofs work in the analogous way in the case where the flux is convex
 and the   corresponding \emph{microscopic convexity} is assumed.

\subsection{Two examples that satisfy microscopic concavity }\label{sc:ex}

Presently we have verified all the hypotheses of Theorem \ref{tm:main} for two
classes of processes. 

\subsubsection{The asymmetric simple exclusion process}\label{sc:asep}

 The asymmetric simple exclusion process (ASEP) was the first example described  in Section 
 \ref{sc:genex}.   It has two parameters $0\le p\ne q\le 1$ such that $p+q=1$.
 To be specific let us take $p>q$ so that on average particles prefer to drift to the right.
  The invariant measure \(\mu^\vr\) is the Bernoulli distribution with parameter
  \(0\le\vr\le1\), while \(\wih\mu^\vr\) is concentrated on zero for any \(\vr\). The hydrodynamic flux is strictly concave: \(\Hc(\vr)=(p-q)\vr(1-\vr)\).

The detailed construction of the processes \(y(t)\) and \(z(t)\) needed for
Assumption \ref{def:microc} can be found in \cite{asepsimple}. Here it is in a nutshell. 

Given the background process \((\un\eta(\cdot),\,\uomp(\cdot))\) and the second class particles \(\{X_m(\cdot)\}\) between them, the processes \(y(\cdot)\) and \(z(\cdot)\) 
are nearest-neighbor random walks on the labels $\{m\}$ with rates $p$ and $q$. 
Walk $y(\cdot)$ has bias to the left (rate $p$ to the left, rate $q$ to the right)
  and walk $z(\cdot)$ has bias to the right (rate $p$ to the right, rate $q$ to the left). 
Their jumps are restricted so that jumps between labels $m$ and $m+1$ are permitted
only when $X_m$ and $X_{m+1}$ are adjacent.  The clocks governing these jumps
are coupled so that the ordering \(y\le z\) is preserved.

 Since a second class particle in ASEP is bounded by a rate one Poisson process, \eqref{eq:expq} holds.
 
We gave an earlier proof of Theorem \ref{tm:main} 
  for ASEP  in \cite{se2/3}.  The present general proof evolved from that earlier one.   

\subsubsection{Totally asymmetric   zero range process with jump rates that increase with exponentially decaying slope} 
\label{sc:1zrp}

This class is completely new in the sequence of models for which \(t^{1/3}\)-scaling of current fluctuations have been verified. Models in this class have a richer behavior than either ASEP or the totally asymmetric zero range process (TAZRP) with constant rate.
As explained in Section \ref{sc:genex},  in a TAZRP
one particle is moved from site $i$ to site $i+1$ at rate $f(\om_i)$, and no
particle jumps to the left (our convention for total asymmetry is \(p=1-q=1\)).  
The jump rate $f:\Zb^+\to\Rb^+$ is nondecreasing,
 $f(0)=0$, and 
$f(z)>0$ for $z>0$. Assume further that $f$ is concave. 

As we shall see later in Section \ref{sc:zrp},
one aspect of microscopic concavity, namely the
ordering of second class  particles, can be achieved
for any TAZRP with a nondecreasing concave jump rate.
Indeed, up to Lemma \ref{lm:zryz} in Section \ref{sc:zrp} we only use monotonicity and concavity of the rates \(f\).
Thus for concave TAZRP only the tail control \eqref{eq:nuc}--\eqref{eq:yzbds}  
of the label processes remains to be provided.
For this part we currently need a stronger hypothesis, detailed in the next assumption. 

\begin{ass}
\label{ass:zrpass}  Let \(p=1-q=1\). Assume the jump rate function $f$ 
of a totally asymmetric zero range process has these properties: 
\begin{itemize}
\item \(f(0)=0<f(1)\),
\item \(f\) is nondecreasing: \(f(z+1)\ge f(z)\),
\item \(f\) is concave with an exponentially decreasing slope: there is an \(0<r<1\) such that for each \(z\ge1\) such that  \(f(z)-f(z-1)>0\),
\be
\frac{f(z+1)-f(z)}{f(z)-f(z-1)}\le r.\label{eq:cvxratio}
\ee
The case where \(f\) becomes constant above some \(z_0\) is included.
\end{itemize}
\end{ass}

\begin{tm} Under  Assumption \ref{ass:zrpass}, a stationary totally asymmetric zero range process
satisfies the conclusions of Theorem \ref{tm:main}, and the conclusions of Corollaries 
\ref{cr:varj}, \ref{cr:lln}, \ref{cr:initial} and \ref{cr:CLT}. 
\label{tm:zrp}\end{tm} 

A class of examples of rates that satisfy Assumption \ref{ass:zrpass} are 
\[
f(z)=1-\exp(-\beta z^\vartheta),\quad \beta>0,\ \vartheta\ge 1.
\]
Another example is the most basic, constant rate TAZRP with  \(f(z)={\bf1}\{z>0\}\).
For this last case a proof has already been given in \cite{julizrp}.

To prove Theorem \ref{tm:zrp} we need to check Assumptions \ref{def:microc} and 
\ref{ass:main} of Theorem \ref{tm:main}.  
The construction of the label processes 
\(y(t)\) and \(z(t)\) and verification of Assumption \ref{def:microc}
 are done in Section \ref{sc:zrp}.
Assumption \ref{ass:main} requires only a few comments. The properties of
the rates required in the first bullet of Assumption \ref{ass:main} are straightforward. 
Since \(f\) is concave and cannot be linear due to \eqref{eq:cvxratio}, Proposition \ref{pr:zrcvx} in the appendix
implies that \(\Hc''(\vr)<0\) for each \(\vr>0\).
Concavity of \(f\) implies bounded jump rates for the second class particle \(Q(t)\), hence a simple Poisson bound gives \eqref{eq:expq}.

\medskip

The remainder of the paper is devoted to proofs. The next two sections prove Theorem
\ref{tm:main}, after that we prove the Strong Law for the second class particle,
and then we return to finish the proof of Theorem \ref{tm:zrp}. 

\section{Upper bound of the main theorem}\label{sc:ub}

In this section we prove the upper bound of \eqref{eq:main}.
We first give a sketch of the proof. As in Section \ref{sc:microconc} on microscopic concavity, we consider the second class particle $Q(t)$ in a pair of processes $(\uomm,\,\uomp)$ at density $\vr$. Additionally there is a positive density of other second class particles that arise from a coupling of $(\uomm,\,\uomp)$ with a third process $\un\eta$ at density $\lambda\in(\vr-\ga_0,\,\vr)$. We emphasize that the coupled $(\un\eta,\,\un\om)$ is not stationary. This is not only because we modified the marginals at the origin to $\wih\mu^\vr$ and $\wih\mu^\la$ from \eqref{eq:muhat}, but more fundamentally because i.i.d.~product measures
%$(\un\eta,\,\un\om)$
are not  stationary for the coupled  evolution. Nevertheless, the marginal processes $\un\om$ and $\un\eta$ are close enough to their stationary  product distributions so that we can calculate conveniently.

The $\om-\eta$ second class particles are conserved during the evolution, and their current is the difference between the currents (heights) of the $\un\om$ and the $\un\eta$ processes. Careful coupling makes it possible to compare $Q(t)$ with the position of a tagged $\om-\eta$ second class particle $X_0(t)$. Fluctuation bounds for $Q(t)$ are derived through several steps: a deviation of $Q(t)$ implies a similar deviation for  $X_0(t)$, which results in a deviation of height differences $h^\om-h^\eta$. The probability of this is bounded by Chebyshev's inequality which brings in variances of the currents $h^\zeta$ and $h^\eta$. These variances are further turned into the first moment of $Q(t)$ essentially via \eqref{eq:varcovar} and \eqref{eq:qchar}. Now the loop is closed, as deviations of $Q(t)$ are bounded by the centered first absolute moment of $Q(t)$. Along the way we see that the sharpest bound is obtained with $\lambda=\vr-c\cdot u/t$ for a constant $c$. We also mention in advance that the critical part of our estimate comes from the order of magnitude $u\sim t^{2/3}$, thus $\vr-\la\sim t^{-1/3}$. With this choice the means for currents and   second class particle velocities that we use for centering provide factors of just the right order for successful completion of the estimation.

 Density $\vr$ is fixed.  Let $\la\in(\vr,\,\vr-\ga_0)$  and apply Assumption \ref{def:microc}
with  constant sequences  \(\vr_i\equiv\vr\) and \(\la_i\equiv\la\) for all \(i\in\Zb\).  Notations 
\(\Pv,\ \Ev,\ \Vv,\ \Cov\) will refer to the coupled four-process evolution
described in Assumption \ref{def:microc}, while
\(\Pv^\vr,\ \Ev^\vr,\ \Vv^\vr,\ \Cov^\vr\) will refer to   a density \(\vr\)
stationary process. Abbreviate 
\be\Psi(t):\,=\Ev|Q(t)-\fl{V^\vr t}|. \label{eq:defPsi}\ee
The requirement that $(\uomm,\,\uomp)$ obey the
basic coupling was included in
Assumption \ref{def:microc}. Consequently $\Psi(t)$ is the $m=1$ expectation of 
 \eqref{eq:main}. 
 
The following lemma does the main work towards the upper bound.
We keep \(\Hc''(\vr)\) explicitly in the estimates, because its non-vanishing is the key feature behind the \(t^{1/3}\)-fluctuations.
\begin{lm}
There exist positive constants \(\al_1,\,\al_2,\,t_0\) such that for each \(t>t_0\)
and   integer $u$
such that 
\(\al_2\sqrt t<u<\al_1t\),
\be
\Pv\{Q(t)>\fl{V^\vr t}+u\}\le C_5\frac{t^2\Hc''(\vr)^2}{u^4}\bigl\{\Psi(t)+u\bigr\}+C_4\frac{t^2}{u^3}.
\label{eq:UB1}\ee
\label{lm:UB1}\end{lm}
\begin{proof}
We start with an integer \(u>0\), and write
\be
\Pv\{Q(t)>\fl{V^\vr t}+u\}\le\Pv\{y(t)\ge k\}+\Pv\{X_k(t)\ge Q(t)>\fl{V^\vr t}+u\}.\label{eq:startw}
\ee
The event \(\{X_k(t)>\fl{V^\vr t}+u\}\) implies that among the $X_m$'s  at most particles
\(X_1,\,\dots,\,X_{k-1}\) have passed the path
\(\bigl(s(\fl{V^\vr t}+u)+1/2\bigr)_{0\le s\le1}\) from right to left.
Each such passing decreases \(h^{\omp}_{\fl{V^\vr t}+u}(t)-h^{\eta}_{\fl{V^\vr t}+u}(t)\) by one (recall 
the statement around  \eqref{eq:2hno}). Hence
we can bound the probability in \eqref{eq:startw} by
\[
\Pv\{y(t)\ge k\}+\Pv\{h^{\omp}_{\fl{V^\vr t}+u}(t)-h^{\eta}_{\fl{V^\vr t}+u}(t)>-k\}.
\]

We   introduce two more processes: \(\uetae\) is a stationary 
process started with initial data \(\etae_i(0)=\eta_i(0)\) for \(i\ne0\),
while \(\etae_0(0)\) is \(\mu^\la\) distributed independently of everything.
\(\uome\) is a stationary 
process started with \(\ome_i(0)=\omp_i(0)\) for \(i\ne0\),
and \(\ome_0(0)\) is \(\mu^\vr\) distributed independently of everything.
Include these   in the basic coupling of $(\un\eta,\,\uomp)$ and write 
\[
\ba
h^{\omp}_{\fl{V^\vr t}+u}(t)-h^{\eta}_{\fl{V^\vr t}+u}(t)
&=h^{\ome}_{\fl{V^\vr t}+u}(t)-h^{\etae}_{\fl{V^\vr t}+u}(t)\\
&\quad+h^{\omp}_{\fl{V^\vr t}+u}(t)-h^{\ome}_{\fl{V^\vr t}+u}(t)\\
&\quad-h^{\eta}_{\fl{V^\vr t}+u}(t)+h^{\etae}_{\fl{V^\vr t}+u}(t).
\ea
\]
Basic coupling implies
\[
\ba
h^{\omp}_{\fl{V^\vr t}+u}(t)-h^{\ome}_{\fl{V^\vr t}+u}(t)
&\le|\omp_0(0)-\ome_0(0)|\le|\omp_0(0)|+|\ome_0(0)|\\
\text{ and} \quad h^{\etae}_{\fl{V^\vr t}+u}(t)-h^{\eta}_{\fl{V^\vr t}+u}(t)
&\le|\etae_0(0)-\eta_0(0)|\le|\etae_0(0)|+|\eta_0(0)|.
\ea
\]
We bound the stationary expectations using  \eqref{eq:Eh}, \eqref{eq:chardef} and Taylor's formula.
This is a crucial computation, which shows that on the characteristic position (that would be case $u=0$), expectation of the height difference is only $\mathcal O(\vr-\lambda)^2$, without constant and first-order expression of the densities.
\begin{align*}
&\quad\ \Ev^\vr h^{\ome}_{\fl{V^\vr t}+u}(t)-\Ev^\la h^{\etae}_{\fl{V^\vr t}+u}(t)\\
&=\Hc(\vr)t-(\fl{V^\vr t}+u)\vr-\Hc(\la)t+(\fl{V^\vr t}+u)\la\\
&\le t\bigl(\Hc(\vr)-\Hc(\la)+\Hc'(\vr)(\la-\vr)\bigr)+u(\la-\vr)+C_1\\
&\le-\frac t2\Hc''(\vr)(\vr-\la)^2+u(\la-\vr)+C_2t(\vr-\la)^3+C_1.
\end{align*}
 $\Hc$ can be differentiated arbitrarily many times,
as we show in Section \ref{sc:afx} of the Appendix. Constant  \(C_1\) above
bounds errors from discarded integer parts. Recall that
tilde stands for the centered random variable.
Collecting terms we continue from \eqref{eq:startw} as follows.
\begin{multline*}
\Pv\{Q(t)>\fl{V^\vr t}+u\}\\
\ba
&\le\Pv\{y(t)\ge k\}\\
&\quad+\Pv\{{\wt h}^{\ome}_{\fl{V^\vr t}+u}(t)-{\wt h}^{\etae}_{\fl{V^\vr t}+u}(t)
>-k+\frac t2\Hc''(\vr)(\vr-\la)^2+u(\vr-\la)\\
&\qquad\qquad-C_2t(\vr-\la)^3-C_1-|\eta_0(0)|-|\etae_0(0)|-|\om_0(0)|-|\ome_0(0)|\}\\
&\le\Pv\{y(t)\ge k\}\\
&\quad+\Pv\{{\wt h}^{\ome}_{\fl{V^\vr t}+u}(t)-{\wt h}^{\etae}_{\fl{V^\vr t}+u}(t)
>\frac t2\Hc''(\vr)(\vr-\la)^2+\frac u2(\vr-\la)\}\\
&\quad+\Pv\{|\eta_0(0)|+|\etae_0(0)|+|\om_0(0)|+|\ome_0(0)|\\
&\qquad\qquad>-k+\frac u2(\vr-\la)-C_2t(\vr-\la)^3-C_1\}.
\ea
\end{multline*}
From now on we use the specific assumption \(\Hc''(\vr)<0\). We maximize the terms
on the right-hand side of the probability of \(\wt h\)'s by the choice
\[
\vr-\la=\frac{-u}{2t\Hc''(\vr)}.
\]
To stay within the range of densities covered
by Assumption \ref{def:microc} we must ensure that \(\la>\vr-\ga_0\).  
So we  introduce a small constant $\al_1>0$ 
and  restrict our calculations
to the case \(u<\al_1 t\). Then
\[
\ba
\Pv\{Q(t)>\fl{V^\vr t}+u\}&\le\Pv\{y(t)\ge k\}\\
&\quad+\Pv\{{\wt h}^{\ome}_{\fl{V^\vr t}+u}(t)-{\wt h}^{\etae}_{\fl{V^\vr t}+u}(t)
>\frac{-u^2}{8t\Hc''(\vr)}\}\\
&\quad+\Pv\{|\eta_0(0)|+|\etae_0(0)|+|\om_0(0)|+|\ome_0(0)|\\
&\qquad>-k-\frac{1}{4\Hc''(\vr)}\cdot\frac{u^2} t+\frac{C_2}{\Hc''(\vr)^3}\cdot\frac{u^3}{t^2}-C_1\}.
\ea
\]
Now we set
\[
k=\fl{\frac{-1}{8\Hc''(\vr)}\cdot\frac{u^2} t},
\]
and assume \(\al_2\sqrt t<u<\al_1 t\) for a possibly smaller \(\al_1\) and a large
enough  \(\al_2\).
That allows us to unify the right-hand side of the inequality in the last line.  Thus 
for all large \(u\) and \(t\) with \(\al_2\sqrt t<u<\al_1 t\)
\[
\ba
\Pv\{Q(t)>\fl{V^\vr t}+u\}&\le\Pv\{y(t)\ge\fl{\frac{-1}{8\Hc''(\vr)}\cdot\frac{u^2} t}\}\\
&\quad+\Pv\{{\wt h}^{\ome}_{\fl{V^\vr t}+u}(t)-{\wt h}^{\etae}_{\fl{V^\vr t}+u}(t)
>\frac{-u^2}{8t\Hc''(\vr)}\}\\
&\ +\Pv\{|\eta_0(0)|+|\etae_0(0)|+|\om_0(0)|+|\ome_0(0)|>C_3\frac{u^2} t\}.
\ea
\]

Assumption \ref{eq:nuc} allows us
to bound the first probability on the right by \(C_4t^2/u^3\) (take \(\ga=\vr-\la\)).
Apply Chebyshev's inequality on the second line and
Markov's inequality on the third one:
\begin{multline*}
\Pv\{Q(t)>\fl{V^\vr t}+u\}\\
\ba
&\le64\frac{t^2\Hc''(\vr)^2}{u^4}\Vv(h^{\ome}_{\fl{V^\vr t}+u}(t)
-h^{\etae}_{\fl{V^\vr t}+u}(t))+C_3\frac t{u^2}+C_4\frac{t^2}{u^3}\\
&\le128\frac{t^2\Hc''(\vr)^2}{u^4}\Bigl\{\Vv^\vr(h^{\ome}_{\fl{V^\vr t}+u}(t))
+\Vv^\la(h^{\etae}_{\fl{V^\vr t}+u}(t))\Bigr\}\\
&\quad+C_4\frac{t^2}{u^3}.
\ea
\end{multline*}
The term \(C_3t/u^2\) was subsumed under  \(C_4t^2/u^3\) due to the condition \(u<\al_1t\).
The variances here are taken under the stationary distributions
of the processes \(\uetae\) and \(\uome\). That allows us 
to apply \eqref{eq:varcovar}, whose right-hand side  takes us back
to the four-process coupling under measure $\Pv$. Recall \eqref{eq:defPsi}. 
\begin{multline*}
\Pv\{Q(t)>\fl{V^\vr t}+u\}\\
\ba
&\le C_5\frac{t^2\Hc''(\vr)^2}{u^4}\Bigl\{\Ev|Q(t)-\fl{V^\vr t}-u|+\Ev|Q^\eta(t)-\fl{V^\vr t}-u|\Bigr\}+C_4\frac{t^2}{u^3}\\
&\le C_5\frac{t^2\Hc''(\vr)^2}{u^4}\Bigl\{\Ev|Q(t)-\fl{V^\vr t}|+\Ev|Q^\eta(t)-\fl{V^\vr t}|+2u\Bigr\}+C_4\frac{t^2}{u^3}\\
&= C_5\frac{t^2\Hc''(\vr)^2}{u^4}\Bigl\{\Psi(t)+2u+\Ev|Q^\eta(t)-\fl{V^\vr t}|\Bigr\}+C_4\frac{t^2}{u^3}.
\ea
\end{multline*}
The variable \(Q^\eta(t)\) above is the location of a single
discrepancy between the process \(\un\eta\) and one started
initially with \(\uetap(0)=\un\eta(0)+\un\de_0\). 

It remains to relate $\Ev|Q^\eta(t)-\fl{V^\vr t}|$ to $\Psi(t)$. This is where 
part \eqref{eq:qq} of  Assumption \ref{def:microc} is a key point. 
Compute now in the four-process coupling 
  of \(\un\eta,\,\uetap,\,\uomm,\,\uomp\) described in Assumption \ref{def:microc}. 
Use \eqref{eq:qq} and Taylor expansion of $\Hc$ again: 
\begin{align}
\Ev|Q^\eta(t)-\fl{V^\vr t}|&\le\Ev(Q^\eta(t)-Q(t))+\Psi(t)\notag\\
&=(\Hc'(\la)-\Hc'(\vr))t+\Psi(t)\label{eq:eqabs}\\
&\le\Hc''(\vr)\cdot(\la-\vr)t+C_6(\vr-\la)^2t+\Psi(t)\notag\\
&=\frac u2+C_6\frac{u^2}t+\Psi(t)\le (\tfrac14+C_6\al_1)u +\Psi(t). \notag
\end{align}
The last inequality used \(u< \al_1 t\). 
Substitute this back into the previous display and rename constants. 
This finishes the proof of \eqref{eq:UB1} and completes the Lemma. 
\end{proof}
Completely analogous arguments lead to the same upper bound for the lower tail of 
$Q(t)$, and together we get the following bound on the tail of the absolute deviation, still for
\(\al_2\sqrt t<u<\al_1t\):
\[
\Pv\{|Q(t)-\fl{V^\vr t}|>u\}\le C_5\frac{t^2\Hc''(\vr)^2}{u^4}\bigl\{\Psi(t)+u\bigr\}+C_4\frac{t^2}{u^3}.
\]

Next we relax the restriction to integral  \(u\) and the upper limit on it:
\begin{lm}
There are positive constants \(\al_2,\,t_0\) such that for all \(t>t_0\) and all \emph{real} \(u>\al_2\sqrt t\),
\[
\Pv\{|Q(t)-\fl{V^\vr t}|>u\}\le C_5\frac{t^2\Hc''(\vr)^2}{u^4}\bigl\{\Psi(t)+u\bigr\}
+C_4\frac{t^2}{u^3}.
\]
\end{lm}
\begin{proof}
Any \(u\ge1\) is less than twice its integer part.
Hence by simply increasing the constants $C_i$,  for all large \(t\) and all \emph{real} \(u\in(\al_2\sqrt t,\,\al_1 t)\),
\be
\Pv\{|Q(t)-\fl{V^\vr t}|>u\}\le C_5\frac{t^2\Hc''(\vr)^2}{u^4}\bigl\{\Psi(t)+u\bigr\}+C_4\frac{t^2}{u^3}.\label{eq:qdev}
\ee

Recall \eqref{eq:expq}. When \(\al_1<\al_0+2|V^\vr|+2\), assume \(\al_1t\le u<(\al_0+2|V^\vr|+2)t\).
Then \(\al_2\sqrt{t}<u\cdot\al_1/(\al_0+2|V^\vr|+2)<\al_1t\) for large enough \(t\),
and \eqref{eq:qdev} still holds for \(u\) replaced by \(u\cdot\al_1/(\al_0+2|V^\vr|+2)\):
\[
\ba
\Pv\{|Q(t)-\fl{V^\vr t}|>u\}&\le\Pv\Bigl\{|Q(t)-\fl{V^\vr t}|>u\cdot\frac{\al_1}{\al_0+2|V^\vr|+2}\Bigr\}\\
&\le C_5\frac{t^2\Hc''(\vr)^2}{u^4}\bigl\{\Psi(t)+u\bigr\}+C_4\frac{t^2}{u^3}
\ea
\]
via modifying the constants by factors of \(\al_1/(\al_0+2|V^\vr|+2)\).

Finally, the case \(u\ge(\al_0+2|V^\vr|+2)t\) will not be relevant for us hence, due to the fact that \(u-|\fl{V^\vr t}|>\al_0t\), we can use \eqref{eq:expq}:
\be\begin{split}
\Pv\{|Q(t)-\fl{V^\vr t}|>u\}&\le\Pv\{|Q(t)|>u-|\fl{V^\vr t}|\}\\
&\le C_7\frac{t^2}{(u-|\fl{V^\vr t}|)^3}\le C_8\frac{t^2}{u^3}.
\end{split}\ee

Combining the above cases  we get the statement for all \(u>\al_2\sqrt t\).
\end{proof}
\begin{proof}[Proof of the upper bound of Theorem \ref{tm:main}]
We now fix \(r>0\), \(1\le m<3\), and write
\begin{multline*}
\Ev\bigl(|Q(t)-\fl{V^\vr t}|^m\bigr)\\
\ba
&=\int_0^\infty\Pv\{|Q(t)-\fl{V^\vr t}|^m>v\}\di v\\
&\le r^mt^{\frac23m}+m\int_{rt^{2/3}}^\infty\Bigl(C_5\frac{t^2\Hc''(\vr)^2}{u^4}\bigl\{\Psi(t)+u\bigr\}+C_4\frac{t^2}{u^3}
\Bigr)u^{m-1}\di u\\
&=r^mt^{\frac23m}+\frac{m C_5\Hc''(\vr)^2}{4-m}r^{m-4}t^{\frac23m-\frac23}\Psi(t)+\frac{m C_5\Hc''(\vr)^2+C_4}{3-m}r^{m-3}t^{\frac23m}.
\ea
\end{multline*}
First choose \(m=1\) and \(r\) large enough to get 
 \(\Psi(t)\le Ct^{2/3}\). Then insert this bound  back into the last line of the display 
 to get the bound for general \(1\le m<3\).
\end{proof}

\def\aaeta{\xi}
\def\abeta{\eta}
\section{Lower bound of the main theorem}\label{sc:lb}

We begin again with an informal preview of the proof. The proof of the lower bound of \eqref{eq:main} uses similar ideas as the upper bound proof, but with an extra twist. The starting point is a pair of processes $(\un\xi,\,\un\xi^+)$ at density $\lambda$ with one second class particle $Q^{(-n)}(t)$ between them started from position $-n$. Coupled to this pair is a process $\un\zeta\ge\un\xi$ that is mostly in density $\vr>\lambda$, except that we set $\un\zeta=\un\xi$ on the interval $-n+1,\,-n+2,\,\dots,\,0$.
%, where its initial density is artificially decreased from $\vr$ to $\lambda$, and so $\zeta=\xi$ on this interval. 
The position $-n$ is chosen so that the $\lambda$-characteristic $-n+V^\lambda t$ started from $-n$ satisfies $V^\lambda t-n=V^\vr t-u$ for a large enough $u>0$ so that the upper bound makes the event $Q^{(-n)}(t)<V^\vr t$ likely.  
%  is smaller by some quantity $u$ than the $\vr$-characteristics $V^\vr t$ that starts from the origin. Now we concentrate on the single $\xi^+-\xi$ second class particle that starts at $-n$. Either it is to the left or to the right of $V^\vr t$ at time $t$, and the aim is to show that both cases lead to deviations of types we need for proving the lower bound. In fact, since our second class particle is expected to be at position $V^\lambda t-n=V^\vr t-u$, we simply bound the probability that the second class particle is on the right of $V^\vr t$ by bounds similar to those we used in the upper bound proof. 
Reasoning as we did for the upper bound, from this event we can deduce an inequality for the current difference between the $\un\zeta$ and the $\un\xi$ processes. 
%The argument that proves this is again similar to those seen in the upper bound. 
In order to turn this inequality into a deviation that can be bounded by Chebyshev's 
inequality as in the upper bound proof, we change the $\un\zeta$ process into a stationary process by introducing the appropriate Radon-Nikodym density for the initial distribution. As in the upper bound proof, the useful perturbation of density is of the order $\vr-\lambda=bt^{-1/3}$. 
% Should we consider}
%$\un\zeta$
%\hl{in a stationary Bernoulli($\vr$) distribution, such bound would not mean a deviation, and the lower bound would not follow. This is where the artificial change of}
%$\un\zeta$\hl{'s}
%\hl{density from $\vr$ to $\lambda$ over sites $-n+1,\,-n+2,\,\dots,\,0$ plays a role. With this modification the bound we obtain becomes in a sense a deviation in current differences, the probability of which is bounded by techniques similar to the upper bound. Finally we show by a Radon-Nikodym type argument that the artificial change in density still allows transforming the estimates to on an unchanged, constant density $\vr$ process. We keep the now-familiar difference $\vr-\lambda=bt^{-1/3}$, and factors that arise from expectations provide the correct scaling of the lower bound.}

Density $\vr$ is fixed again,  and $\la\in(\vr-\ga_0,\,\vr)$ is a varying auxiliary density.
We let the jointly defined four processes $(\un{\abeta},\,\un{\abeta}^+,\,\uomm,\,\uomp)$
be exactly as defined in the upper bound proof of
 Section \ref{sc:ub}, namely, as given by 
  Assumption \ref{def:microc}  with constant  densities $\la_i\equiv\la$ and $\vr_i\equiv\vr$.
 The initial distribution of $(\un{\abeta},\,\uomp)$ is 
${\un{\wih\mu}}^{\un\la,\un\vr}$ of \eqref{eq:unmu}.
Two second class particles start from the origin:
 $Q^{\abeta}$   between processes
    $\un{\abeta}$ and $\un{\abeta}^+$,  and 
\(Q\)  between  processes \(\uomm\) and \(\uomp\). 
The quantity of primary interest is abbreviated, as
before, by $\Psi(t)=\Ev|Q(t)-\fl{V^\vr t}|$.

To prove the lower bound of \eqref{eq:main} it suffices, by Jensen's inequality, 
to prove the case $m=1$. This means showing that
$\Psi(t)\ge Ct^{2/3}$ for large $t$ and a constant $C>0$.

\subsection{Perturbing a segment initially}

For this proof we need to introduce another coupled system and invoke
Assumption \ref{def:microc} once more. 
 By concavity of the flux   characteristic speeds $V^\vr=\Hc'(\vr)$ and
$V^\la=\Hc'(\la)$ satisfy 
  $V^\vr\le V^\la$. Throughout this section  $u>0$ denotes a
positive integer, and
\[
n=\fl{V^\la t}-\fl{V^\vr t}+u.
\]

Recall definitions \eqref{eq:cmuhat} and \eqref{eq:cmu} of the 
single-site coupling measures.  Let 
$(\un\aaeta(\cdot),\,\un\ze(\cdot))$ be a pair of processes that obeys the basic
coupling, and whose initial distribution is the product measure 
\[
\Bigl(\bigotimes_{i<-n}\mu^{\la,\vr}\Bigr)\Bigl(\bigotimes_{i=-n}{\wih\mu}^{\la,\vr}\Bigr)\Bigl(\bigotimes_{-n<i\le0}\mu^{\la,\la}\Bigr)\Bigl(\bigotimes_{0<i}\mu^{\la,\vr}\Bigr).
\]
This initial measure complies with the pattern in \eqref{eq:unmu},
but translated  \(n\) sites to the left so that \({\wih\mu}^{\la,\vr}\) 
is the distribution at site \(-n\) instead of the origin.  A few points about this initial
state:   $\un\aaeta(0)$ has the stationary density-$\lambda$ product distribution except
at site $-n$ where it is $\wih\mu^\la$-distributed.
$\un\ze(0)$ has the product distribution with marginals
$\mu^\vr$, except at sites $\{-n+1,\,\dots,\,0\}$
where the parameter $\vr$ switches to $\la$, and at site
$-n$ where it has distribution $\wih\mu^\vr+1$. 
At sites  \(-n<i\le0\) \(\mu^{\la,\la}\) forces \(\aaeta_i(0)=\ze_i(0)\). 

We add a second class particle to the process $\un \aaeta(\cdot)$, 
 start it at site $-n$ and denote its position at time $t$
by $Q^{(-n)}(t)$. Let  $\un \aaeta^+(t):= \un\aaeta(t)+
\un \delta_{Q^{-n}(t)}$.

As described in Section \ref{sc:microconc} 
 the   $\ze-\aaeta$ second class particles are labeled and their
 ordered positions denoted by $\{X_m(t)\}$.  The 
 labeling is chosen to satisfy initially 
\be
\dots\le X_{-1}(0)\le X_0(0)=-n<0<X_1(0)\le X_2(0)\le\dots
\label{eq:LBinitX}\ee
Thus initially
$X_0(0)=-n=Q^{(-n)}(0)$. We invoke Assumption \ref{def:microc} to have a label
process $z(t)$ with tail bound \eqref{eq:yzbds} such that 
$Q^{(-n)}(t)=X_{z(t)}(t)$. (Here $\un\aaeta$ plays the role of $\un\eta$ and 
 $\un\zeta$ plays the role of $\uomp$ of Assumption
\ref{def:microc}).  

As before, the heights (or currents, recall \eqref{eq:hno}) of the processes
$\un\aaeta(\cdot)$ and $\un\ze(\cdot)$ are denoted by
$h^{\aaeta}_{\fl{Vt}}$ and $h^{\ze}_{\fl{Vt}}$, respectively.   The first observation is that
 $Q^{(-n)}$ gives one-sided control over the difference
of these currents.

\begin{lm}\label{lm:geom}
For any $i\in\Zb$
\[
Q^{(-n)}(t)\leq i \quad\text{implies}\quad
h^{\ze}_i(t)-h^{\aaeta}_i(t)\leq -z(t).
\]
\end{lm}
\begin{proof}
Recall
again, from \eqref{eq:hno} and the statement around \eqref{eq:2hno},  that the 
height difference \(h^{\ze}_i(t)-h^{\aaeta}_i(t)\)  equals the net
number of second class  particle  passings  of the path \(\bigl(si+1/2\bigr)_{0\le s\le1}\)
from left to right.  That is, each left-to-right passing increases 
 \(h^{\ze}_i(t)-h^{\aaeta}_i(t)\) while each right-to-left passing decreases it.

 Suppose \(z(t)\le0\). Then  \eqref{eq:LBinitX} and 
\(X_{z(t)}(t)=Q^{(-n)}(t)\leq i\) imply that only those second class particles with labels 
\(z(t)+1,\,z(t)+2,\,\dots,\,0\) could have passed the path \(\bigl(si+1/2\bigr)_{0\le s\le1}\)
from left to right.  The claim follows. 
 
If \(z(t)>0\), then \(X_{z(t)}(t)=Q^{(-n)}(t)\leq i\) implies that at least 
those second class particles with labels 
\(1,\,2,\,\dots,\,z(t)\) have crossed the path
\(\bigl(si+1/2\bigr)_{0\le s\le1}\)
from right to left.  Again the claim follows. 
\end{proof}

Let $\un{\wih\om}(\cdot)$ be a process started from the product distribution
\(\bigl(\,\bigotimes\limits_{i\ne-n}\mu^\vr\bigr)\otimes(\wih\mu^\vr+1)\).
The next lemma compares the initial distributions of $\un\ze$ and
$\un{\wih\om}$.  No coupling of $\un\ze$ and
$\un{\wih\om}$ is proposed or required. 
\begin{lm}\label{lm:rn}
There exist constants  $\ga=\gamma(\vr)>0$ 
and $C_1(\vr)<\infty$ 
such that for all
$\la\in(\vr-\gamma,\,\vr)$ and all events $A$ the following inequality holds:
\[
\Pv\{\un\ze\in A\}\leq\Pv\{ \un{\wih\om}\in A\}^\frac12\cdot\exp\bigl\{C_1(\vr)
n(\vr-\la)^2 \bigr\}.
\]
\end{lm}
\begin{proof}

We use the Cauchy-Schwarz inequality below to perform a change of
measure on the distribution of the $\un\zeta$ process. First we condition on
the initial \(\un\ze\)-configuration at sites $\{-n+1,\,-n+2, \dots, -1,\,0\}$.
\[
\ba &\Pv\{\un\ze\in A\}=\sum_{z_{-n+1},\dots z_{-1},
z_0}\Pv\{\,\un\ze\in A\,|\,\zeta_{-n+1}(0)=z_{-n+1},\dots,
\zeta_{0}(0)=z_0\}\\
&\qquad\qquad\times\Big[\prod_{i=-n+1}^0 \mu^\vr(z_i)\Big]^\frac12 
 \prod_{i=-n+1}^0\frac{\mu^\la(z_i)}{[\mu^\vr(z_i)]^\frac12}\\
&\leq\Bigl[\sum_{z_{-n+1},\dots,z_0}[\Pv\{\,\un\ze\in A\,|\,\zeta_{-n+1}(0)=z_{-n+1},\dots,
\zeta_{0}(0)=z_0\}]^2\,
\prod_{i=-n+1}^0\mu^\vr(z_i)\Bigr]^\frac12\\
&\qquad\qquad\times\Bigl[\sum_{z_{-n+1},\dots,z_0}
\prod_{i=-n+1}^0 \frac{[\mu^\la(z_i)]^2}
{\mu^\vr(z_i)}\Bigr]^\frac12\\
&\leq\Bigl[\sum_{z_{-n+1},\dots,z_0}\Pv\{\,\un\ze\in A\,|\,\zeta_{-n+1}(0)=z_{-n+1},\dots,
\zeta_{0}(0)=z_0\}\,
\prod_{i=-n+1}^0\mu^\vr(z_i)\Bigr]^\frac12\\
&\qquad\qquad\times\Bigl[\sum_{z_{-n+1},\dots,z_0}
\prod_{i=-n+1}^0 \frac{[\mu^\la(z_i)]^2}
{\mu^\vr(z_i)}\Bigr]^\frac12\\
&=\Pv\{ \un{\wih\om}\in A\}^\frac12\cdot\Bigl[\sum_{z_{-n+1},\dots,z_0}
\prod_{i=-n+1}^0 \frac{[\mu^\la(z_i)]^2}{\mu^\vr(z_i)}\Bigr]^\frac12.
\ea
\]
The last inequality came from dropping the square. 
For the last equality note that the distributions of the initial configurations $\{\wih\om_i(0)\}$ and 
$\{\zeta_i(0)\}$ are product-form and agree outside the interval $\{-n+1,\,-n+2, \dots, -1,\,0\}$.
Thus conditioned on the initial values in $\{-n+1,\,-n+2, \dots, -1,\,0\}$ these processes
have identical conditional probabilities.  

To complete the proof we bound the last factor
in brackets. 
 Recall formulas \eqref{eq:zdef} and \eqref{eq:mudef} for the state sum and the site-marginals. 
 Without the power $1/2$ the factor in brackets equals 
\[\sum_{z_{-n+1},\dots z_0}
\Bigl(\frac{Z(\te(\vr))}{Z(\te(\la))^2}\Bigr)^n \prod_{i=-n+1}^0
\frac{\e{(2\te(\la)-\te(\vr))
z_i}}{f(z_i)!}=\Bigl(\frac{Z(2\te(\la)-\te(\vr))Z(\te(\vr))}{Z(\te(\la))^2}\Bigr)^{n}.\]
In the appendix we show that $\log Z(\te)$ and $\te(\vr)$ are infinitely differentiable.  
Let  $\ve=\te(\vr)-\te(\la)$. By local Lipschitz continuity of the function \(\te(\vr)\), the interval \((\te(\la)-\ve,\,\te(\la)+\ve)\) is in \((\un\te,\,\bar\te)\) with a small enough choice of \(\ga\). There exists some \(\te\in(\te(\la)-\ve,\,\te(\la)+\ve)\) such that 
\[
\ba
\log\Big(\frac{Z(2\te(\la)-\te(\vr))Z(\te(\vr))}{Z(\te(\la))^2}\Big)&=
\log Z(\te(\la)-\ve)+\log Z(\te(\la)+\ve)\\
&\quad-2\log Z(\te(\la))\\
&=\frac12\frac{\di^2}{\di
\te^2}\log Z(\te) \ve^2
\le C_1(\vr)\cdot(\vr-\la)^2.
\ea
\]
Thus we get the bound 
\[\Big(\frac{Z(2\te(\la)-\te(\vr))Z(\te(\vr))}{Z(\te(\la))^2}\Big)^n
\le \exp\{C_1(\vr)\cdot n(\vr-\la)^2\}.
\qedhere
\]
\end{proof}

\subsection{Completion of the proof of the lower bound}
The gist of the proof is to get upper bounds on the complementary probabilities
$\Pv\{Q^{(-n)}(t)>\fl{V^\vr t}\}$ and $\Pv\{Q^{(-n)}(t)\le \fl{V^\vr t}\}$.
 As stated
$u$ is an arbitrary but positive integer
 and $n=\fl{V^\la t}-\fl{V^\vr t}+u$.
\begin{lm}
\[
\Pv\{Q^{(-n)}(t)>\fl{V^\vr
t}\}\leq\frac{\Psi(t)}{u}+\frac{C_2t(\vr-\la)}{u}+\frac{2}{u}.
\]
\label{lm:LBlm3} \end{lm}
\begin{proof}
Distributionwise the system $(\un\aaeta,\,\un\aaeta^+,\,Q^{(-n)})$ is a translate of
   $(\un{\abeta},\,\un{\abeta}^+,\,Q^{\abeta}) $, and so 
\begin{align*}
&\Pv\{Q^{(-n)}(t)>\fl{V^\vr t}\}
=\Pv\{Q^{(-n)}(t)+n-\fl{V^\la t}>u\}\\&\qquad =\Pv\{Q^{\abeta}(t)-\fl{V^\la
t}>u\}  \leq\frac{\Ev(|Q^{\abeta}(t)-\fl{V^\la t}|)}{u}\\
&\qquad  \leq\frac{\Ev(|Q^{\abeta}(t)-Q(t)|)}{u}
+\frac{\Ev(|Q(t)-\fl{V^\vr
t}|)}{u}+\frac{\fl{V^\la t}-\fl{V^\vr t}}{u}.  \end{align*}
Use \eqref{eq:qq} precisely as was done in \eqref{eq:eqabs}
to conclude that the first term equals
\[ u^{-1}\Ev(Q^{\abeta}(t)-Q(t)) = u^{-1}{t}\big(H'(\la)-H'(\vr)\big)=-{u}^{-1}H''(\nu)t(\vr-\la)\]
for some $\nu\in(\la,\,\vr)$. The second term is $\Psi(t)/u$, and the
third term is similarly estimated by
$-u^{-1}{H''(\nu)}t(\vr-\la)+{2}/{u}$, the last part coming
from discarded integer parts. Setting
$C_2:=2\max\limits_{\nu\in[\vr-\gamma,\,\vr]}-H''(\nu)$ finishes the proof.
\end{proof}
Notice that \(H''(\vr)<0\) was crucial in the previous proof, as well as in the following lemma, and the final proof thereafter. These points show where the proof fails for symmetric systems -- recall that these would have lower-order current fluctuations on the characteristics.
\begin{lm} Let $K=K(t)$ satisfy $0<K<-\frac13tH''(\vr)(\vr-\la)^2$.
Then for small enough $\ga>0$,  
  large enough  $t$, and $\la\in(\vr-\ga,\,\vr)$, 
\[
\ba \Pv\{Q^{(-n)}(t)\leq\fl{V^\vr t}\}&\leq
\frac{\Vv^\vr(\om_0)^{1/2}\Psi(t)^{1/2}}{-\frac13tH''(\vr)(\vr-\la)^2-K}
\cdot\e{C_1 n(\vr-\la)^2}\\
&\quad+\frac{C_4}{-\frac16
tH''(\vr)(\vr-\la)^2 -C_3t(\vr-\la)^3-\vr}\cdot\e{C_1 n(\vr-\la)^2} \\
&\quad+ \frac{\Vv^\la(\eta_0)\Psi(t)}{K^2/4}+\frac{C_6 t
(\vr-\la)}{K^2} +\frac{C_5}{K-4|\la|}\\
&\quad+Ct^{\kappa-1}\ga^{2\kappa-3}K^{-\kappa}.
\ea
\]
\label{lm:LBlm4} \end{lm}
\begin{proof}
Lemma \ref{lm:geom} leads to
\begin{align}
\Pv\{Q^{(-n)}(t)\leq\fl{V^\vr
t}\}&\leq\Pv\{h^{\ze}_{\fl{V^\vr t}}(t)-h^{\aaeta}_{\fl{V^\vr t}}(t)
\leq -z(t)\}\notag\\
&\leq\Pv\{-z(t)\ge K/4\}\label{eq:lab}\\
&\quad+\Pv\{h^{\ze}_{\fl{V^\vr t}}(t)\leq K+t\big(H(\la)-\la H'(\vr)\big)\}\label{eq:torn}\\
&\quad+\Pv\Bigl\{h^{\aaeta}_{\fl{V^\vr t}}(t)>3K/4+t\big(H(\la)-\la
H'(\vr)\big)\Bigr\}. \label{eq:tovarj}
\end{align}
To bound \eqref{eq:lab} we use the assumed distribution bound \eqref{eq:nuc}
on $z(t)$ and get
\[
\Pv\{-z(t)\ge K/4\} \le Ct^{\kappa-1}\ga^{2\kappa-3}K^{-\kappa}.
\]
Apply Lemma \ref{lm:rn} to line \eqref{eq:torn} to
bound it by the pro\-ba\-bi\-li\-ty of the process $\un{\wih\om}$:
\[
(\ref{eq:torn})\le \bigl[\Pv\{h^{\wih\om}_{\fl{V^\vr t}}(t)\leq
K+t\big(H(\la)-\la H'(\vr)\big)\}\bigr]^\frac{1}{2}
\cdot \e{C_1 n(\vr-\la)^2}.
\]
As in the proof of Lemma \ref{lm:UB1} we compare with a coupled stationary processes to get
precise bounds:
\[
\ba
h^{\wih\om}_{\fl{V^\vr t}}(t)&=\wt h^{\ome}_{\fl{V^\vr t}}(t)+[h^{\wih\om}_{\fl{V^\vr t}}(t)-h^{\ome}_{\fl{V^\vr t}}(t)]\\
&\quad+[\Ev h^{\ome}_{\fl{V^\vr t}}(t)- t(H(\vr)-\vr H'(\vr))]+ t(H(\vr)-\vr H'(\vr))\\
&\ge\wt h^{\ome}_{\fl{V^\vr t}}(t)-|\wih\om_{-n}(0)|-|\ome_{-n}(0)|-|\vr|+t(H(\vr)-\vr H'(\vr)).
\ea
\]
After the equality sign, the absolute value of the first
 term in brackets is not larger  than $|\wih\om_{-n}(0)-\ome_{-n}(0)|\le
|\wih\om_{-n}(0)|+|\ome_{-n}(0)|$.   The second term in brackets  is between
$-|\vr|$ and $|\vr|$ due to the integer part in \(\fl{V^\vr t}\).
Consequently 
\[
h^{\wih\om}_{\fl{V^\vr t}}(t)\leq K+t(H(\la)-\la H'(\vr))
\]
implies 
\[
\ba
\wt h^{\ome}_{\fl{V^\vr t}}(t) -|\wih\om_{-n}(0)|-|\ome_{-n}(0)|&\le K+t[H(\la)-H(\vr)+H'(\vr)(\vr-\la)]+|\vr|\\
&\leq K+\frac12tH''(\vr)(\vr-\la)^2+ C_3t(\vr-\la)^3+|\vr|.
\ea
\]
Then, we cut the
event into two parts  according to the value of 
 $|\wih\om_{-n}(0)|+|\ome_{-n}(0)|$ and we use \eqref{eq:varcovar} to bound the variance of 
$\Vv [h^{\ome}_{\fl{V^\vr t}}(t)]$ by the
function $\Psi(t)$.
\[
\ba
(\ref{eq:torn})&\leq\bigl[\Pv^{\vr}\{\wt h^{\ome}_{\fl{V^\vr t}}(t)\leq K+\frac13
tH''(\vr)(\vr-\la)^2\}\bigr]^\frac12
\cdot\e{C_1 n(\vr-\la)^2}\\
&\quad+
\bigl[\Pv\{|\wih\om_{-n}(0)|+|\ome_{-n}(0)|>-\frac{1}{6}tH''(\vr)(\vr-\la)^2-C_3t(\vr-\la)^3-|\vr|\}\bigl]^\frac12\\
&\quad\quad\cdot\e{C_1
n(\vr-\la)^2}\\
&\leq\frac{\Vv^\vr(h^{\ome}_{\fl{V^\vr t}}(t))^{1/2}}{-\frac13tH''(\vr)(\vr-\la)^2-K}
\cdot\e{C_1 n(\vr-\la)^2}\\
&\quad+\frac{[\Ev(|\wih\om_{-n}(0)|+|\ome_{-n}(0)|)^2]^{1/2}}{-\frac16
tH''(\vr)(\vr-\la)^2 -C_3t(\vr-\la)^3-|\vr|}\cdot\e{C_1 n(\vr-\la)^2}\\
&\leq\frac{\Vv^\vr(\om_0)^{1/2}\Psi(t)^{1/2}}{-\frac13tH''(\vr)(\vr-\la)^2-K}
\cdot\e{C_1 n(\vr-\la)^2}\\
&\quad+\frac{C_4}{-\frac16
tH''(\vr)(\vr-\la)^2 -C_3t(\vr-\la)^3-|\vr|}\cdot\e{C_1 n(\vr-\la)^2}.
\ea
\]
Now we turn to \eqref{eq:tovarj}.
To reduce $h^{\aaeta}_{\fl{V^\vr t}}$ to the current of the density-$\la$ 
equilibrium process $h^{\etae}_{\fl{V^\vr t}}$ and to get rid of the integer
part errors we argue as before. 
\[
\ba
h^{\aaeta}_{\fl{V^\vr t}}&=\wt h^{\etae}_{\fl{V^\vr t}}+[h^{\aaeta}_{\fl{V^\vr t}}-h^{\etae}_{\fl{V^\vr t}}]\\
&\quad+[\Ev^\la h^{\etae}_{\fl{V^\vr t}}-t(H(\la)-\la H'(\vr))]+ t(H(\la)-\la H'(\vr)).
\ea
\]
$h^{\aaeta}_{\fl{V^\vr t}}(t)$ differs by
at most $|\aaeta_{-n}(0)-\etae_{-n}(0)|\le |\aaeta_{-n}(0)|+|\etae_{-n}(0)|$ from \linebreak[4] $h^{\etae}_{\fl{V^\vr t}}(t))$.  Taking integer parts again into account,
giving another error term $|\la|$, line \eqref{eq:tovarj} is bounded
from above by
\[
\Pv\Bigl\{\wt h^{\etae}_{\fl{V^\vr t}}(t)+|\aaeta_{-n}(0)|+|\etae_{-n}(0)|+|\la|
\ge 3K/4\Bigr\}.
\]
 Then, we cut the event into two parts and use Markov's inequality
 in the second one:
\[
\ba
\eqref{eq:tovarj}&\le\Pv^\la\Bigl\{\wt h^{\etae}_{\fl{V^\vr t}}(t)\ge K/2\Bigr\}+
\Pv\Bigl\{|\aaeta_{-n}(0)|+|\etae_{-n}(0)|> K/4-|\la|\Bigr\}\\
&\leq
\frac{\Vv^\la(h^{\etae}_{\fl{V^\vr t}})}{K^2/4}+\frac{C_5}{K-4|\la|}.
\ea
\]
We can use \eqref{eq:varcovar} again to continue with
\[
\eqref{eq:tovarj}\le\frac{\Vv^\la(\aaeta_0)\Ev(|Q^{\abeta}(t)-\fl{V^\vr t}|)}{K^2/4}+\frac{C_5}{K-4|\la|}.
\]
Repeating the first two steps of calculation \eqref{eq:eqabs} we can write 
\[
\ba
\Ev(|Q^{\abeta}(t)-\fl{V^\vr t}|)&\le\Ev(|Q^{\abeta}(t)-Q(t)|)+\Ev(|Q(t)-\fl{V^\vr t}|)\\
&\le Ct(\vr-\la)+\Psi(t). 
\ea
\]
So, we finally get
\[ \eqref{eq:tovarj}\le \frac{\Vv^\la(\eta_0)\Psi(t)}{K^2/4}+\frac{C_6 t (\vr-\la)}{K^2}  +\frac{C_5}{K-4|\la|}.
\qedhere
\]
\end{proof}

\begin{proof}[Proof of the lower bound of Theorem \ref{tm:main}]
As observed in the beginning of this Section, 
  it suffices to prove  that
\be
\liminf_{t\to\infty} t^{-2/3}\Psi(t)>0.
\label{eq:LB8}\ee
In the last two lemmas take
\[
u=\lceil ht^{2/3}\rceil,\quad\vr-\la=bt^{-1/3}, \quad \text{ and
}\quad K=bt^{1/3},
\]
where $h$ and  $b$ are large, in particular $b$ large enough to have
$b<-\frac13 H''(\vr) b^2$ so that $K$ satisfies the assumption of
Lemma \ref{lm:LBlm4}. Then
\[
\ba n&=\fl{V^\la t}-\fl{V^\vr t}+u\\
&\le(H'(\la)-H'(\vr))t+u+2\\
&=-H''(\vr)(\vr-\la)t+ u+ C_7 t (\vr-\la)^2+2\\
&\le(-H''(\vr)b+h) t^{2/3} + C_7 b^2 t^\frac13 +3\\
&\le C_8 t^{2/3} \ea
\] for large enough $t$. With these definitions we can simplify the outcomes of Lemma
\ref{lm:LBlm3} and Lemma \ref{lm:LBlm4}  to the inequalities \be
\Pv\{Q^{(-n)}(t)>\fl{V^\vr t}\}\leq C\frac{\Psi(t)}{t^{2/3}}
+\frac{C_2b}{h} +\frac{2}{h t^{2/3}} \label{eq:LBaux3}\ee and
\be\begin{split} \Pv\{Q^{(-n)}(t)\leq\fl{V^\vr t}\}&\leq
C\biggl(\frac{\Psi(t)}{t^{2/3}}\biggr)^{1/2}
+C\frac{\Psi(t)}{t^{2/3}}+ \frac{C_6}{b}+ \frac{C_5}{b t^\frac13}+
Cb^{\kappa-3}.
\end{split}
\label{eq:LBaux4}\ee The new constant $C$ depends on $b$ and $h$.

The lower bound \eqref{eq:LB8}  now follows because the left-hand sides of
\eqref{eq:LBaux3}--\eqref{eq:LBaux4} add up to $1$ for each fixed
$t$, while we can fix
 $b$ large enough and then $h$ large enough so
that $C_2b/h + C_6/b+Cb^{\kappa-3} < 1$ (recall \(\kappa<3\)).  Then $t^{-2/3}\Psi(t)$ must have a
positive lower bound for all large enough $t$. This completes the
proof of Theorem \ref{tm:main}.
\end{proof}

\section{Strong Law of Large Numbers for the second class particle}
\label{sc:Qlln}
This section proves  the Strong Law of Large Numbers (Corollary \ref{cr:lln}). We assume that the
 jump rates of the second class particle are bounded, i.e.,
\be
\left.\ba
&p(y+1,\,z)-p(y,\,z),&p(y,\,z)-p(y,\,z+1)\\
&q(y,\,z+1)-q(y,\,z),&q(y,\,z)-q(y+1,\,z)
\ea\right\}
\le C \quad \forall \omin\le y,\,z<\omax.\label{eq:boundrate}
\ee
This means that the
second class particle has at most rate $C$ to jump to the right and to
the left, respectively, implying that starting at any time \(t\), it can be bounded by rate $C$ Poisson processes that start from its position \(Q(t)\).
\begin{proof}[Proof of Corollary \ref{cr:lln}]
 Let $\ve,\,\de>0$. Define the events \[A_n:=\Bigl\{
\Bigl|\frac{Q(n^{1+\de})}{n^{1+\de}}-V^\vr\Bigr|>\ve/2\Bigr\}\]
for \(n\in\Nb\).
Then,
Markov's inequality and Theorem \ref{tm:main} imply, for \(1\le m<3\) and large \(n\),
\[ \ba \Pv\{A_n\} &= \Pv \{ \big|Q(n^{1+\de})-V^\vr n^{1+\de}\big|^m >
(\ve/2)^m n^{(1+\de)m}\}\\
&\le \frac{1}{(\ve/2)^m n^{(1+\de)m}}\cdot\Ev[|Q(n^{1+\de})-V^\vr
n^{1+\de}|^m] \\
&\le \frac{C_1}{(3-m)(\ve/2)^m}\cdot\frac{1}{n^{m(1+ \de)/3}}, \ea
\] which is summable if $(1+\de)m>3$. Here
$\de$ can be chosen arbitrarily small by taking \(m\) close to 3. By the Borel-Cantelli Lemma
 there exists  a.s.\ $n_0\in \mathbb N$
such that \be \label{n1+d} \forall n\ge n_0 \qquad
\Big|\frac{Q(n^{1+\de})}{n^{1+\de}}-V^\vr\Big|<\ve/2.\ee 

Using this we show that a.s.\ there exists $n_1\in \mathbb N$ such
that  \be
\label{eq:slln}\Big|\frac{Q(t)}{t}-V^\vr\Big|<\ve
\quad\text{for all real $t\ge n_1^{(1+\de)}$.}
\ee 

Let $n\ge n_0$ and suppose there exists
some $t \in [n^{1+\de},\,(n+1)^{1+\de})$ such that
\eqref{eq:slln} fails: $|Q(t)-V^\vr t|\ge \ve t$.
Together with  \eqref{n1+d} we have, if $n$ is large,  
\be
\ba
|Q(t)-Q(n^{1+\de})|
&\ge|Q(t)-V^\vr t|-|Q(n^{1+\de})-V^\vr n^{1+\de}|-|V^\vr t-V^\vr n^{1+\de}|\\
&\ge\ve t-\ve/2\cdot n^{1+\de}-|V^\vr|(t-n^{1+\de})\\
&\ge\frac{\ve}{4}n^{1+\de}. 
\ea\label{eq:bad}
\ee

The jump rates \eqref{eq:boundrate} (both left and right) of \(Q\) are bounded by \(C\). However, the event \eqref{eq:bad} implies that at least \(\fl{\frac{\ve}{4}n^{1+\de}}\) many left jumps or this many right jumps happen in the time interval \([n^{1+\de},\,(n+1)^{1+\de})\). For large \(n\), the length of this interval is smaller than \(2(1+\de)n^\de\). Let \(N(\cdot)\) be a rate \(C\) Poisson process. Then for large \(n\) the probability of the event \eqref{eq:bad} is bounded from above by
\[
\ba
2\Pv\{N(2(1+\de) n^\de)\ge\frac{\ve}{4}n^{1+\de}\}
&\le2\Pv\{e^{N(2(1+\de)n^\de)}\ge\e{\ve/4\cdot n^{1+\de}}\}\\
&\le2\e{-\ve/4\cdot n^{1+\de}}\Ev[e^{N(2(1+\de)n^\de)}]\\
&=2\e{-\ve/4\cdot n^{1+\de}}\cdot\e{(\e{}-1)2C(1+\de)n^{\de}}.
\ea
\]
This quantity is summable over $n$, so the Borel Cantelli Lemma implies that 
a.s.\ \eqref{eq:slln} holds eventually.
Since this is true for each \(\ve>0\),   the Strong Law of Large Numbers holds.
\end{proof}

\section{Microscopic concavity for a class of totally asymmetric concave exponential zero range processes}\label{sc:zrp}

In this section we verify that Assumption \ref{def:microc} can be satisfied under
Assumption  \ref{ass:zrpass}, and thereby complete  the proof of Theorem \ref{tm:zrp}.
 
 The task is to construct   the processes \(y(t)\) and \(z(t)\) with the requisite properties.
 First let the processes \((\un\eta(\cdot),\,\uomp(\cdot))\) 
 evolve in the basic coupling so that \(\eta_i(t)\le\omp_i(t)\) for all $i\in\Zb$ and $t\ge 0$.
  We consider as a background process this pair with the labeled and ordered  $\omp-\eta$ 
 second class particles \(\dotsm\le X_{-2}(t)\le X_{-1}(t)\le X_0(t)\le X_1(t)\le X_2(t)\le\dotsm\).
  
  At each time \(t\ge0\) this background induces a partition \(\{\Mc_i(t)\}\) of the label
  space \(\Zb\)
  into intervals  indexed by sites \(i\in\Zb\), with partition intervals given by 
\[
\Mc_i(t):\,=\{m\,:\,X_m(t)=i\}.
\]
(For simplicity we assumed infinitely many second class particles in both directions, but no problem arises in case we only have finitely many of them.)   \(\Mc_i(t)\) contains the labels of 
the second class particles that reside at site \(i\) at time \(t\), and can be empty. 
The labels of the second class particles that are at the same site as the one labeled \(m\) 
form the set \(\Mc_{X_m(t)}(t)=\,:\{a^m(t),\,a^m(t)+1,\,\dots,\,b^m(t)\}\). 
The processes $a^m(t)$ and $b^m(t)$ are always well-defined and satisfy
$a^m(t)\le m\le b^m(t)$. 

Let us clarify these notions by discussing the ways in which \(a^m(t)\) and \(b^m(t)\) can change.
\begin{itemize}
\item A second class particle jumps from site \(X_m(t-)-1\) to site \(X_m(t-)\). Then this one necessarily has label \(a^m(t-)-1\), and it becomes the lowest labeled one at site \(X_m(t-)=X_m(t)\) after the jump. Hence \(a^m(t)=a^m(t-)-1\).
\item A second class particle, different from \(X_m\), jumps from site \(X_m(t-)\) to site \(X_m(t-)+1\). Then this one is necessarily labeled \(b^m(t-)\), and it leaves the site \(X_m(t-)\), hence \(b^m(t)=b^m(t-)-1\).
\item The second class particle \(X_m\) is the highest labeled on its site, that is, \(m=b^m(t-)\), and it jumps to site \(X_m(t-)+1\). Then this particle becomes the lowest labeled in the set \(\Mc_{X_m(t-)+1}=\Mc_{X_m(t)}\), hence \(a^m(t)=m\). In this case \(b^m(t)\) can be computed from \(b^m(t)-a^m(t)+1=\omp_{X_m(t)}(t)-\eta_{X_m(t)}(t)\), the number of second class particles at the site of \(X_m\) after the jump.
\end{itemize}

We fix initially \(y(0)=z(0)=0\).  The evolution of $(y,\,z)$ is superimposed on
the background evolution $(\un\eta,\,\uomp,\,\{X_m\})$   following the general rule below: 
Immediately after every move of the background process that involves the site where
$y$ resides before this move,   $y$ picks a new value from the labels on the site where it 
resides after the move.   
Thus $y$ itself jumps only within partition
intervals $\Mc_i$. But $y$  joins a  new partition interval   whenever it is the
highest $X$-label on its site and its 
  ``carrier''
particle $X_y$ is  forced to move to the next site on the right.  This is the situation when
$y(t-)=b^{y(t-)}(t-)$ and at time  $t$ an $\omp-\eta$ move from this site
happens. 
(Recall that the choice of $X$-particle to move is determined by rule \eqref{eq:Xrule}. In
the present case there is only one type of $\omp-\eta$ move: the highest label from a site moves
to the next site on the right.)  All this works for $z$ in exactly the same way.  

Next we specify the probabilities that $y$ and $z$ use to refresh their values.
When $y$ and $z$ reside at separate sites, they refresh independently.  When they
are together in the same partition interval, they use the joint distribution 
 in the third bullet below. 

\begin{itemize}
\item Whenever any change occurs in either \(\uomp\) or \(\un\eta\) at site \(X_{y(t-)}(t-)\) and, as a result of the jump,
\(a^{y(t-)}(t)\ne a^{z(t-)}(t)\), that is, \(y(t-)\) and \(z(t-)\) belong to different parts \emph{after} the jump then, independently of everything else,
\be
y(t):\,=\left\{\ba
&a^{y(t-)}(t),&&\text{with pr.\ }\frac{f(\omp_{X_{y(t-)}(t)}(t)-1)-f(\eta_{X_{y(t-)}(t)}(t))}{f(\omp_{X_{y(t-)}(t)}(t))-f(\eta_{X_{y(t-)}(t)}(t))},\\
&b^{y(t-)}(t),&&\text{with pr.\ }\frac{f(\omp_{X_{y(t-)}(t)}(t))-f(\omp_{X_{y(t-)}(t)}(t)-1)}{f(\omp_{X_{y(t-)}(t)}(t))-f(\eta_{X_{y(t-)}(t)}(t))}
\ea\right.\label{eq:ych}
\ee
when the denominator is non-zero, and \(y(t):\,=a^{y(t-)}(t)\) when the denominator is zero.
\item Whenever any change occurs in either \(\uomp\) or \(\un\eta\) at site \(X_{z(t-)}(t-)\) and, as a result of the jump,
\(a^{y(t-)}(t)\ne a^{z(t-)}(t)\), that is, \(y(t-)\) and \(z(t-)\) belong to different parts \emph{after} the jump then, independently of everything else,
\be
z(t):\,=\left\{\ba
&b^{z(t-)}(t)-1,&&\text{with pr.\ }\frac{f(\omp_{X_{z(t-)}(t)}(t))-f(\eta_{X_{z(t-)}(t)}(t)+1)}{f(\omp_{X_{z(t-)}(t)}(t))-f(\eta_{X_{z(t-)}(t)}(t))},\\
&b^{z(t-)}(t),&&\text{with pr.\ }\frac{f(\eta_{X_{z(t-)}(t)}(t)+1)-f(\eta_{X_{z(t-)}(t)}(t))}{f(\omp_{X_{z(t-)}(t)}(t))-f(\eta_{X_{z(t-)}(t)}(t))}
\ea\right.\label{eq:zch}
\ee
when the denominator is non-zero, and \(z(t):\,=b^{z(t-)}(t)\) when the denominator is zero. When \(\omp_{X_{z(t-)}(t)}(t)=\eta_{X_{z(t-)}(t)}(t)+1\),\quad\(b^{z(t-)}(t)-1\) is not an admissible value but in this case the probability in the first line is zero.
\item Whenever any change occurs in either \(\uomp\) or \(\un\eta\) at sites \(X_{y(t-)}(t-)\) or \(X_{z(t-)}(t-)\) and, as a result of the jump,
\(a^{y(t-)}(t)=a^{z(t-)}(t)\), that is, \(y(t-)\) and \(z(t-)\) belong to the same part \emph{after} the jump, that is, \(X_{y(t-)}(t)=X_{z(t-)}(t)\) then, independently of everything else,
\be
\begin{pmatrix}y(t)\\z(t)\end{pmatrix}:\,=\left\{\ba
&\begin{pmatrix}a^{y(t-)}(t)\\b^{y(t-)}(t)-1\end{pmatrix},\\
&\quad\text{with pr.\ }\frac{f(\omp_{X_{y(t-)}(t)}(t))-f(\eta_{X_{y(t-)}(t)}(t)+1)}{f(\omp_{X_{y(t-)}(t)}(t))-f(\eta_{X_{y(t-)}(t)}(t))},\\
&\begin{pmatrix}a^{y(t-)}(t)\\b^{y(t-)}(t)\end{pmatrix},\\
&\quad\text{with pr.\ }\frac{f(\eta_{X_{y(t-)}(t)}(t)+1)-f(\eta_{X_{y(t-)}(t)}(t))}{f(\omp_{X_{y(t-)}(t)}(t))-f(\eta_{X_{y(t-)}(t)}(t))}\\
&\qquad\text{\phantom{with pr.\ }}-\frac{f(\omp_{X_{y(t-)}(t)}(t))-f(\omp_{X_{y(t-)}(t)}(t)-1)}{f(\omp_{X_{y(t-)}(t)}(t))-f(\eta_{X_{y(t-)}(t)}(t))},\\
&\begin{pmatrix}b^{y(t-)}(t)\\b^{y(t-)}(t)\end{pmatrix},\\
&\quad\text{with pr.\ }\frac{f(\omp_{X_{y(t-)}(t)}(t))-f(\omp_{X_{y(t-)}(t)}(t)-1)}{f(\omp_{X_{y(t-)}(t)}(t))-f(\eta_{X_{y(t-)}(t)}(t))}
\ea\right.\label{eq:yzch}
\ee
when the denominator is non-zero, and
\[
(y(t),\,z(t)):\,=(a^{y(t-)}(t),\,b^{y(t-)}(t))
\]
when the denominator is zero. When \(\omp_{X_{z(t-)}(t)}(t)=\eta_{X_{z(t-)}(t)}(t)+1\),\quad\(b^{z(t-)}(t)-1\) is not an admissible value but in this case the probability in the first line is zero.
\end{itemize}
The fact that the numbers on the right hand-sides are probabilities follows from \(\omp_i(t)>\eta_i(t)\) on the sites \(i\) in question, and from the monotonicity and concavity of \(f\). 
The above moves for \(y\) and \(z\) always occur within labels   at a given site. 
This determines whether the particle
 \(Q(t):\,=X_{y(t)}(t)\) or  \(Q^\eta(t):\,=X_{z(t)}(t)\)  is the one to jump if the next move 
 out of the site is   an $\omp-\eta$ move. 
 
We prove that the above construction has the properties required in
Assumption \ref{def:microc}.
\begin{lm}
The pair $(\uomm,\,\uomp):\,=(\uomp-\un\delta_{X_y},\,\uomp)$ obeys basic coupling, 
as does the pair $(\un\eta,\,\uetap):\,=(\un\eta,\,\un\eta+\un\delta_{X_z})$. 
\end{lm}
\begin{proof}
We write the proof for \((\uomm,\,\uomp)\). We need to show that, given the 
configuration \((\un\eta,\,\uomp,\,\{X_m\},\,y)\), the jump rates of \((\uomm,\,\uomp)\) 
are the ones prescribed in  basic coupling (Section \ref{sc:basiccoupling}) and by 
\(\eqref{eq:add}\).  
Leftward jumps of type \eqref{eq:rem} do not happen in the system under discussion. 
Since the jump rate function $p$ depends only on its first argument, jumps out of sites 
$i\ne Q$ 
 happen for \(\uomm\) and \(\uomp\) with the same rate 
 \(p(\om^-_i,\,\om^-_{i+1})=f(\om^-_i)=f(\om_i)=p(\om_i,\,\om_{i+1})\). 
  The only point to consider  is jumps out of site \(i=Q\).

Since the last time any change occurred at site \(i\), \(y\) chose values according to \eqref{eq:ych} or \eqref{eq:yzch}. Notice that \eqref{eq:ych} and \eqref{eq:yzch} give the same marginal probabilities for this choice. Hence 
\be
\text{\(y\) took on value \(a^y\) with probability}\quad 
\frac{f(\omp_i-1)-f(\eta_i)}{f(\omp_i)-f(\eta_i)}\label{eq:ynjpr}
\ee
and
\be \text{\(y\) took on value \(b^y\) with probability}\quad 
\frac{f(\omp_i)-f(\omp_i-1)}{f(\omp_i)-f(\eta_i)},\label{eq:yjpr}
\ee
as given in \eqref{eq:ych}, or $y$ took on value  \(a^y\) in the case  \(f(\omp_i)=f(\eta_i)\). 
According to the basic coupling of \(\un\eta\) and \(\uomp\), the following jumps can occur over the edge \((i,\,i+1)\):
\begin{itemize}
\item With rate \(p(\omp_i,\,\omp_{i+1})-p(\eta_i,\,\eta_{i+1})=f(\omp_i)-f(\eta_i)\), when positive, \(\uomp\) jumps without \(\un\eta\). The highest labeled second class particle, \(X_{b^y}\) jumps from site \(i\) to site \(i+1\).
\begin{itemize}
\item With probability \eqref{eq:yjpr} \(X_y=Q\) jumps with \(X_{b^y}\). In this case
\[
\omm_i(t-)=\omp_i(t-)-1=\omp_i(t)=\omm_i(t)
\]
since the difference \(Q\) disappears from site \(i\). Also,
\[
\omm_{i+1}(t-)=\omp_{i+1}(t-)=\omp_{i+1}(t)-1=\omm_{i+1}(t),
\]
since the difference \(Q\) appears at site \(i+1\). So in this case 
$\uomp$ undergoes a jump but    \(\uomm\) does not, and the rate is 
\[
[f(\omp_i)-f(\eta_i)]\cdot\frac{f(\omp_i)-f(\omp_i-1)}{f(\omp_i)-f(\eta_i)}=f(\omp_i)-f(\omm_i).
\]

\item With probability \eqref{eq:ynjpr} \(X_y=Q\) does not jump with \(X_{b^y}\), 
since it has label \(a^y\) and not \(b^y\) (this probability is zero if \(\omp_i=\eta_i+1\)). In this case \(\uomm\) and  \(\uomp\)
perform the same jump   and it occurs with rate
\[
[f(\omp_i)-f(\eta_i)]\cdot\frac{f(\omp_i-1)-f(\eta_i)}{f(\omp_i)-f(\eta_i)}=f(\omm_i)-f(\eta_i).
\]
\end{itemize}
\item With rate \(p(\eta_i,\,\eta_{i+1})=f(\eta_i)\), both \(\un\eta\) and \(\uomp\) jump over the edge \((i,\,i+1)\). No change occurs in the \(\omp-\eta\) particles, hence no change occurs in \(Q\). This implies that the process \(\uomm\) jumps as well.
\end{itemize}
Summarizing  we see that the rate for \((\uomm,\,\uomp)\) to jump together 
over \((i,\,i+1)\) is \(f(\omm_i)\), and the rate for $\uomp$  to  jump without
$\uomm$ is $f(\omp_i)-f(\omm_i)$.  This is exactly what basic coupling requires. 

A very similar argument can be repeated for \((\un\eta,\,\uetap)\).  
\end{proof}
\begin{lm}\label{lm:zryz}
Inequality \eqref{eq:yz} $y\le z$ holds in the above construction.
\end{lm}
\begin{proof}
Since no jump of $y$ or $z$ moves one of them into a new partition interval, the 
only situation that can jeopardize \eqref{eq:yz} is the simultaneous refreshing 
of $y$ and $z$ in a common partition interval.  But this case is governed by 
step \eqref{eq:yzch}  which by definition ensures that  \(y\le z\). 
\end{proof}
So far in this section everything is valid for a general zero range process with nondecreasing concave jump rate.  Now we use 
the special convexity requirement \eqref{eq:cvxratio}.
With $r\in(0,\,1)$ from \eqref{eq:cvxratio}, define the geometric distribution
\be
\nu(m):\,=\left\{\ba
&(1-r)r^m,&&m\ge0\\
&0,&&m<0.
\ea\right.\label{eq:nugeom}
\ee
 
\begin{lm}\label{lm:zrgeom}
Conditioned on the process \((\un\eta,\,\uomp)\), 
the bounds \(y(t)\overset{\text{d}}{\le}\nu\) and \(z(t)\overset{\text{d}}{\ge}-\nu\) hold 
for all $t\ge 0$.
\end{lm}
The proof of this lemma is achieved in three steps.
\begin{lm}\label{lm:ygeom}
Let \(Y\) be a random variable with distribution \(\nu\), and fix integers  \(a\le b\)
and  \(\eta<\omp\) 
so that \(\omp-\eta=b-a+1\).  Apply the following operation to $Y$: 

{\rm (i)} if \(a\le Y\le b\), apply the probabilities from \eqref{eq:ych}  {\rm (}equivalently,
\eqref{eq:ynjpr} and \eqref{eq:yjpr}{\rm )} with   parameters \(a,\,b,\,\eta,\,\omp\)
 to pick a new value for $Y$; 
 
 {\rm (ii)} if $Y<a$ or $Y>b$ then do not change $Y$. 
 
  Then the resulting distribution \(\nu^*\) is stochastically dominated by \(\nu\).
\end{lm}
\begin{proof}
There is nothing to prove when \(b=a\), hence we assume \(b>a\) or, equivalently, \(\omp-\eta=b-a+1\ge2\). It is also clear that \(\nu^*(m)=\nu(m)\) for \(m<a\) or \(m>b\). We need to prove, in view of the distribution functions,
\[
\sum_{\ell=a}^m\nu^*(\ell)\ge\sum_{\ell=a}^m\nu(\ell)\text{\quad or, equivalently,\quad}\sum_{\ell=m}^b\nu^*(\ell)\le\sum_{\ell=m}^b\nu(\ell)
\]
for all \(a\le m\le b\). Notice that \(\nu^*\) gives zero weight on values \(a<m<b\) (if any), therefore the left hand-side of the second inequality equals \(\nu^*(b)\) for \(a<m\le b\). Hence the above display is proved once we show
\begin{align}
\nu^*(b)&\le\nu(b)\text{, that is,}\notag\\
\frac{f(\omp)-f(\omp-1)}{f(\omp)-f(\eta)}\cdot\sum_{\ell=a}^b\nu(\ell)&\le\nu(b),\label{eq:telefscope}
\end{align}
see \eqref{eq:ych}. When \(f(\omp)=f(\omp-1)\), there is nothing to prove. Hence assume \(f(\omp)>f(\omp-1)\) which by concavity implies that \(f\) has positive increments on \(\{\eta,\,\dots,\,\omp\}\). If \(b<0\) then both sides are zero. If \(b\ge0\) then we have, by \eqref{eq:cvxratio},
\[
\ba
\nu(\ell)\le\nu(b)\cdot r^{\ell-b}&\le\nu(b)\cdot\prod_{z=\omp-b+\ell}^{\omp-1}\frac{f(z)-f(z-1)}{f(z+1)-f(z)}\\
&=\nu(b)\cdot\frac{f(\omp-b+\ell)-f(\omp-b+\ell-1)}{f(\omp)-f(\omp-1)}
\ea
\]
for each \(\ell\le b\). The first inequality also takes into account possible \(\nu(\ell)=0\) values for negative \(\ell\)'s. With this we can write
\[
\sum_{\ell=a}^b\nu(\ell)\le\nu(b)\cdot\frac{f(\omp)-f(\omp-b+a-1)}{f(\omp)-f(\omp-1)}
\]
which becomes \eqref{eq:telefscope} via \(\omp-\eta=b-a+1\).
\end{proof}
We repeat the lemma for \(z(t)\).
\begin{lm}
Let \(Z\) be a random variable of distribution \(-\nu\), and fix integers
\(a\le b\), \(\eta<\omp\)  so that \(\omp-\eta=b-a+1\). Operate on $Z$ as was done for $Y$ in Lemma \ref{lm:ygeom}, but this
time use the probabilities from 
 \eqref{eq:zch} with parameters \(a,\,b,\,\eta,\,\omp\).  Let  \(-\nu^*\) be 
  the resulting distribution.  Then  \(\nu^*\) is stochastically dominated by \(\nu\).
\end{lm}
\begin{proof}
Again, we assume \(b>a\) or, equivalently, \(\omp-\eta=b-a+1\ge2\). It is also clear that \(\nu^*(-m)=\nu(-m)\) for \(m<a\) or \(m>b\). We need to prove
\[
\sum_{\ell=a}^m\nu^*(-\ell)\le\sum_{\ell=a}^m\nu(-\ell)
\]
for all \(a\le m\le b\). Notice that \(-\nu^*\) gives zero weight on values \(a\le\ell<b-1\) (if any), therefore the left hand-side of the inequality equals 0 for \(a\le m<b-1\), \(\nu^*(b-1)\) for \(m=b-1\), and agrees to the right hand-side for \(m=b\). Hence the above display is proved once we show
\begin{align}
\nu^*(-b)&\ge\nu(-b)\text{, that is,}\notag\\
\frac{f(\eta+1)-f(\eta)}{f(\omp)-f(\eta)}\cdot\sum_{\ell=a}^b\nu(-\ell)&\ge\nu(-b),\label{eq:telefscope2}
\end{align}
see \eqref{eq:zch}.
We have, by \eqref{eq:cvxratio},
\[
\ba
\nu(-\ell)\ge\nu(-b)\cdot r^{b-\ell}&\ge\nu(-b)\cdot\prod_{\underset{f(z)>f(z-1)}{z=\eta+1}}^{\eta+b-\ell}\frac{f(z+1)-f(z)}{f(z)-f(z-1)}\\
&=\nu(-b)\cdot\frac{f(\eta+1+b-\ell)-f(\eta+b-\ell)}{f(\eta+1)-f(\eta)}
\ea
\]
for each \(\ell\le b\). The first inequality also takes into account possible \(\nu(-b)=0\) values for positive \(b\)'s. With this we can write
\[
\sum_{\ell=a}^b\nu(-\ell)\ge\nu(-b)\cdot\frac{f(\eta+1+b-a)-f(\eta)}{f(\eta+1)-f(\eta)}
\]
which becomes \eqref{eq:telefscope2} via \(\omp-\eta=b-a+1\).
\end{proof}
\begin{lm}\label{lm:yattr}
The dynamics defined by \eqref{eq:ych} or \eqref{eq:zch} is attractive.
\end{lm}
\begin{proof}
Following the same realizations of \eqref{eq:ych}, we see that two copies of \(y(\cdot)\) under a common environment can be coupled so that whenever they get to the same part \(\Mc_i\), they move together from that moment. The same holds for \(z(\cdot)\).
\end{proof}
\begin{proof}[Proof of Lemma \ref{lm:zrgeom}]
Initially \(y(0)=0\) by definition, which is clearly a distribution dominated by \(\nu\) of \eqref{eq:nugeom}. Now we argue recursively: by time \(t\) the distribution of \(y(t)\) was a.s.\ only influenced by finitely many jumps of the environment, which resulted in distributions \(\nu_1\), then \(\nu_2\), then \(\nu_3\), 
etc. Suppose \(\nu_k\overset{\text{d}}{\le}\nu\), and let \(\nu^*\) be the distribution that would result from \(\nu\) by the \(k+1^\text{st}\) jump. Then \(\nu_{k+1}\overset{\text{d}}{\le}\nu^*\) by \(\nu_k\overset{\text{d}}{\le}\nu\) and Lemma \ref{lm:yattr}, while \(\nu^*\overset{\text{d}}{\le}\nu\) by Lemma \ref{lm:ygeom}. A similar argument proves the lemma for \(z(\cdot)\).
\end{proof}

\appendix
\appendixpage

\section{Convexity and total positivity}
\label{sc:convtotpos}
This section derives a general convexity result for exponentially tilted measures.
Let $\nu$ be a nondegenerate probability measure on $\Rb$ and assume
that for some open interval $\Ic\subseteq\Rb$, 
\be
\Zca(\te)= \int e^{\te x}\,\nu(\dii x)<\infty
\qquad\text{for all $\te\in\Ic$.}  
\label{eq:CAass1}\ee
For $\te\in\Ic$ define the exponentially tilted measures 
$\nu^\te$ by 
\[
\int g\,d\nu^\te = \Zca(\te)^{-1}\int g(x)
e^{\te x}\,\nu(\dii x)
\]
(for bounded Borel test functions $g$).  
The nondegeneracy assumption (that $\nu$ is not supported on a
single point) and  \eqref{eq:CAass1}  guarantee that
\[
\vr(\te)= \int x\,\nu^\te(\dii x)
\]
is a finite, continuous, strictly increasing function that maps $\Ic$ onto 
a nontrivial open interval $\Jc$. For $\vr\in\Jc$  the inverse function
is denoted by $\te(\vr)$.   

Let $\psi$ be a measurable function on $\Rb$, and assume 
(by shrinking $\Ic$ if necessary) that 
\[
\int \lvert{\psi}\rvert\,d\nu^\te<\infty
\qquad\text{for all $\te\in\Ic$.}  
\]
Since $\abs{x}^k\le k!\ve^{-k}(e^{\ve x}+e^{-\ve x})$ for any $\ve>0$
and $\Ic$ is an open interval, it follows
that $\int \lvert{\psi}\rvert \abs{x}^k\,d\nu^\te<\infty$ for all $k\ge 0$
and $\te\in\Ic$.
Consequently  as a  function of $\te$ the integral  
$\int \psi\,d\nu^\te$ has derivatives of all orders.  

A particular case 
is $\psi(x)=x$ which gives the infinite differentiability of $\vr(\te)$. 
Let us also note the infinite differentiability of the inverse function
$\te(\vr)$.  
Since $\vr'(\te)$ is the variance of the distribution $\nu^\te$,
$\vr'(\te)>0$ by the nondegeneracy of $\nu$, and so directly from
the definition of the derivative $\te'(\vr)=1/\vr'(\te(\vr))$. 
Repeated use of basic differentiation rules produces all 
derivatives $\te^{(n)}(\vr)$.
Notice that this argument shows a uniform lower and upper bound of $\vr'(\te)$, that is, Lipschitz continuity of both $\vr(\te)$ and $\te(\vr)$ on bounded closed intervals.

Define
\[
\Psi(\vr)=\int {\psi}\,d\nu^{\te(\vr)}.
\]
 $\Psi$ is also infinitely differentiable as a composite of two such
 functions. 
  
\begin{tm} Assume $\psi$ is a convex function on $\Rb$. Then 
  $\Psi$ 
is convex on $\Jc$.   Assume furthermore that  no
 linear function $g(x)=ax+b$ satisfies 
 $\psi=g$  $\nu$-a.e.  Then  $\Psi''(\vr)>0$ for all $\vr\in\Jc$ 
and in particular $\Psi$ is strictly convex on $\Jc$. 
\label{tm:convtm1}
\end{tm}

\begin{proof}  The proof can be reduced 
to the theory of total positivity.  
In what follows, citations and terminology
 are from Karlin's monograph \cite{totpos}.  The claims made
 in our  Theorem \ref{tm:convtm1} follow
 from applying Theorem 3.5(a)--(c) from p.~285
of \cite{totpos} to the operator 
\[
T\psi(\vr) = \int \psi\di\nu^{\te(\vr)}=\int_{\Rb} K(\vr,\,x)\psi(x)\,\nu(\dii x),
\quad \vr\in\Jc,
\]
where the kernel is defined by  $K(\vr,\,x)=\Zca(\vr(\te))^{-1}e^{\te(\vr)x}$. 
 The  property of the kernel $K$ that gives the result  is {\sl extended total positivity}
(ETP) of   order 3.  This is the requirement of strict positivity on certain types of
determinants of partials of dimensions up to $3\times 3$: for all $(\vr,\,x)\in\Jc\times\Rb$,
\be
K^*\binom{\overset{\text{$n$ entries}}{\overbrace{\vr, \dotsc,\vr}}}{x,\dotsc,x}=
\det_{1\le i,\,j\le n}\,
\Bigl[\frac{\partial^{i+j-2}}{\partial \vr^{i-1}\partial x^{j-1}}
K(\vr,\,x)\Bigr]>0 \quad
 \text{for $n=1,\,2,\,3$.} 
 \label{eq:defETP} \ee
We argue this in stages.

We first observe that the kernel 
$L(\te,\,x)=\Zca(\te)^{-1}e^{\te x}$ on $\Ic\times\Rb$ 
is ETP of all orders.
Recall that the {\sl Wronskian} of $n$ functions $f_1,\dotsc,f_n$ 
is the $n\times n$ determinant 
\[
W[f_1,\dotsc,f_n](x)=\det_{1\le i,\,j\le n}[f_i^{(j-1)}(x)].
\]
If $u$ is another function, the Wronskian satisfies the identity
\be
W[uf_1,\dotsc,uf_n](x)=u(x)^nW[f_1,\dotsc,f_n](x).
\label{eq:Wrident}\ee
To justify \eqref{eq:Wrident}, 
  Leibniz's rule 
\[
(uf_i)^{(j-1)}=\sum_{k=1}^j\binom{j-1}{k-1} f_i^{(k-1)}u^{(j-k)}  \quad (1\le j\le n)
\]
implies that
the matrix $A=[(uf_i)^{(j-1)}(x)]_{1\le i,\,j\le n}$ is the product of the
matrices 
\[B=[(f_i)^{(k-1)}(x)]_{1\le i,\,k\le n}\quad\text{and}\quad  
C=\bigl[ \tbinom{j-1}{k-1}u^{(j-k)}(x)\ind\{k\le j\}\bigr]_{1\le k,\,j\le n}. 
\]
By upper-triangularity $\det C=u(x)^n$. Then 
the corresponding determinant identity    $\det (A)=\det (B)\cdot\det (C)$ is precisely 
\eqref{eq:Wrident}.
  
 Now we can verify the ETP property of kernel $L$, utilizing \eqref{eq:Wrident}:
\begin{align*}
&\det_{1\le i,\,j\le n}\,
\Bigl[\frac{\partial^{i+j-2}}{\partial x^{i-1}\partial \te^{j-1}}
L(\te,\,x)\Bigr]=
\det_{1\le i,\,j\le n}\,
\Bigl[\frac{\partial^{j-1}}{\partial \te^{j-1}}
\bigl\{\Zca(\te)^{-1}\te^{i-1}e^{\te x}\bigr\}\Bigr]\\
&\qquad=\Zca(\te)^{-n}e^{n\te x}  
W[1,\,\te, \dotsc, \te^{n-1}]\\
&\qquad=\Zca(\te)^{-n}e^{n\te x} \prod_{j=1}^{n-1} j! >0.  
\end{align*}   
  
To go from $L(\te,\,x)$ to $K(\vr,\,x)=L(\te(\vr),\,x)$, consider
the $3\times 3$ determinant that appears in \eqref{eq:defETP}, apply the chain rule
and a row operation:
\begin{align*}
&\begin{vmatrix} K&K_x&K_{xx}\\K_{\vr}&K_{\vr x}&K_{\vr xx}\\
K_{\vr\vr}&K_{\vr\vr x}&K_{\vr\vr xx}\\
\end{vmatrix}\\
&\qquad=
\begin{vmatrix} L&L_x&L_{xx}\\L_{\te}\te_\vr&L_{\te x}\te_\vr&L_{\te xx}\te_\vr\\
L_{\te\te}\te_\vr^2+L_{\te}\te_{\vr\vr}
&L_{\te\te x}\te_\vr^2+L_{\te x}\te_{\vr\vr}
&L_{\te\te xx}\te_\vr^2+L_{\te xx}\te_{\vr\vr}\\
\end{vmatrix}\\
&\qquad
=\begin{vmatrix} L&L_x&L_{xx}\\L_{\te}\te_\vr&L_{\te x}\te_\vr&L_{\te xx}\te_\vr\\
L_{\te\te}\te_\vr^2 &L_{\te\te x}\te_\vr^2 
&L_{\te\te xx}\te_\vr^2 \\
\end{vmatrix}
=\te_\vr^3\begin{vmatrix} L&L_x&L_{xx}\\L_{\te}&L_{\te x}
&L_{\te xx}\\
L_{\te\te}&L_{\te\te x}&L_{\te\te xx} \\
\end{vmatrix}>0.
\end{align*}
The last inequality is by the ETP property of kernel $L$ and the strict positivity 
$\te_\vr>0$ of the derivative.  The $1\times 1$ and $2\times 2$
determinants in \eqref{eq:defETP} are principal minors of the determinant above
and are positive by the same reasoning.
  
  We have shown that the kernel $K$ has the ETP property of  order 3.  
In addition to ETP, Theorem 3.5 from p.~285
of \cite{totpos} requires the hypotheses  
\[
\int K(\vr,\,x)\,\nu(\dii x)=1 
\quad\text{and}\quad
\int K(\vr,\,x)x\,\nu(\dii x)= a\vr+b
\]
for some $a>0$ and $b\in\Rb$. The first
one is true  by virtue of the 
normalization $\Zca(\te)^{-1}$,    and the second one 
with $a=1$ and $b=0$ by the definition of $\vr(\te)$. 
The proof is now completed by an appeal to Theorem 3.5 from p.~285
of \cite{totpos}.
\end{proof}

These convexity properties can also be proved in an elementary
way by developing suitable correlation inequalities. 
Such a proof is given in the note \cite{convex}. 
We are indebted to an anonymous referee of that note for
pointing out the connection with total positivity. 

Subsequent sections of the appendix 
  extract from Theorem \ref{tm:convtm1} 
 consequences for the processes we study.

\section{Monotonicity of measures}\label{sc:cvxi}

In this part of the appendix we show that the measures \(\mu^\vr\) and \(\wih\mu^\vr\) defined in \eqref{eq:mudef} and \eqref{eq:muhat}, respectively, are stochastically monotone
as functions of $\vr$. We start with a simple
\begin{lm}\label{lm:dercov}
Fix a function \(\vp(\om)\) on \(\Zb\), 
bounded by a polynomial. Then \newline
\(\Ev^\te(\vp(\om))\) is differentiable in \(\te\) on \((\un\te,\,\bar\te)\), and
\[
\frac{\di}{\di\te}\Ev^\te(\vp(\om))=\Cov^\te(\vp(\om),\,\om).
\]
\end{lm}
\begin{proof}
Convergence of the series involved in \(\Ev^\te(\vp(\om))\) can be verified via the ratio test, even after differentiating the terms. Since \(\mu^\te\) is the exponentially weighted version of \(\mu^{\te_0}\) for some \(\te_0\), we have
\[
\ba
\frac{\di}{\di\te}\,\Ev^\te\vp(\om)&=\frac{\di}{\di\te}\,\frac{\Ev^{\te_0}(\vp(\om)\cdot\e{(\te-\te_0)\om})}{\Ev^{\te_0}\e{(\te-\te_0)\om}}\\
&=\frac{\Ev^{\te_0}(\vp(\om)\cdot\om\cdot\e{(\te-\te_0)\om})}{\Ev^{\te_0}\e{(\te-\te_0)\om}}-\Ev^{\te_0}(\vp(\om)\cdot\e{(\te-\te_0)\om})\cdot\!\frac{\Ev^{\te_0}(\om\cdot\e{(\te-\te_0)\om})}{[\Ev^{\te_0}\e{(\te-\te_0)\om}]^2}\\
&=\Cov^\te(\vp(\om),\,\om).
\ea\qedhere
\]
\end{proof}
\begin{cor}
For any \(\un\te<\te<\bar\te\), the state sum \eqref{eq:zdef} satisfies
\begin{align}
\frac{\di}{\di\te}\log Z(\te)&=\frac{1}{Z(\te)}\sum_{z=\omin}^{\omax}z \frac{e^{\te z}}{f(z)!}=\Ev^\te(\om)=\,:\vr(\te),\label{eq:d1z}\\
\frac{\di^2}{\di\te^2}\log Z(\te)&=\frac{\di}{\di\te}\vr(\te)=\Vv^\te(\om).\label{eq:d2z}
\end{align}
The function \(\vr(\te)\) is strictly increasing   and 
maps  \((\un\te,\,\bar\te)\) onto \((\omin,\,\omax)\).
\end{cor}
\begin{proof} Everything is already covered except the last surjectivity statement. 
Due to the monotonicity and continuity one only needs to show convergence at
the boundaries $\un\te$, $\bar\te$  to $\omin$, $\omax$. First let us consider the case when \(\bar\te<\infty\). Then \(\omax=\infty\) and Fatou's lemma implies
\[
\liminf_{\te\nearrow\bar\te}Z(\te)=\liminf_{\te\nearrow\bar\te}\sum_{z\in I}\frac{\e{\te z}}{f(z)!}\ge\sum_{z\in I}\liminf_{\te\nearrow\bar\te}\frac{\e{\te z}}{f(z)!}=\sum_{z\in I}\frac{\e{\bar\te z}}{f(z)!}=\infty
\]
since for \(z>0\)
\[
\frac{\e{\bar\te z}}{f(z)!}=\prod_{y=1}^z\frac{\e{\bar\te}}{f(y)}\ge1
\]
by definition of \(\bar\te\) and \(f\) being nondecreasing. This shows that \(\log Z(\te)\) takes on arbitrarily large values as \(\te\nearrow\bar\te\). We also know that it is a smooth and convex function on \((\un\te,\,\bar\te)\) (see \eqref{eq:d2z}). This implies that its derivative \eqref{eq:d1z} is not bounded from above i.e., arbitrarily large \(\vr\) values can be achieved. The same reasoning works in case \(\un\te>-\infty\) for arbitrarily large negative \(\vr\) values.

When \(\bar\te=\infty\) then, regardless whether \(\omax\) is finite or infinite, fix any \(0\le y<\omax\) and write
\be
\ba
\vr(\te)&=\Ev^\te(\om\cdot{\bf1}\{\om>y\})+\Ev^\te([\om]^+\cdot{\bf1}\{\om\le y\})-\Ev^\te([\om]^-\cdot{\bf1}\{\om\le y\})\\
&\ge(y+1)\cdot\Pv^\te(\om>y)-\Ev^\te([\om]^-\cdot{\bf1}\{\om\le y\})\\
&\ge(y+1)-(y+1)\cdot\Pv^\te(\om\le y)-\sqrt{\Ev^\te(([\om]^-)^2)}\cdot\sqrt{\Pv^\te(\om\le y)}\\
&\ge(y+1)-(y+1)\cdot\Pv^\te(\om\le y)-\sqrt{\Ev^{\te_0}(([\om]^-)^2)}\cdot\sqrt{\Pv^\te(\om\le y)}
\ea\label{eq:rhosplit}
\ee
for a fixed \(\un\te<\te_0<\te\). The last inequality follows by monotonicity of \(\mu^\te\) in \(\te\) and \(([\om]^-)^2)\) being a nonincreasing function of \(\om\).
For any \(\omin-1<z\le y\) and \(\te>\un\te\),
\[
\frac{\mu^\te(z)}{\mu^\te(y+1)}=\prod_{x=z}^y\frac{\mu^\te(x)}{\mu^\te(x+1)}=\prod_{x=z}^y\frac{f(x+1)}{\e{\te}}\le\Bigl(\frac{f(y+1)}{\e{\te}}\Bigr)^{y-z+1}.
\]
Given \(0\le y<\omax\) and \(1>\ve>0\), there is a large enough \(\te\) which makes the last fraction smaller than \(\ve\). With such a choice we have
\[
\Pv^\te\{\om\le y\}=\sum_{z=\omin}^y\mu^\te(z)\le\mu^\te(y+1)\sum_{z=\omin}^y\ve^{y-z+1}\le\ve\cdot\frac{1-\ve^{y-\omin+1}}{1-\ve}.
\]
Therefore, for the case of a finite \(\omax\), choosing \(y=\omax-1\) and large \(\te\) makes \eqref{eq:rhosplit} arbitrarily close to \(\omax\). When \(\omax=\infty\), the argument shows that \(\vr(\te)\ge y+1\) can be achieved for any \(y\ge0\). A similar computation demonstrates that any density towards \(\omin\) can be reached when \(\un\te=-\infty\).
\end{proof}

\begin{cor}\label{cr:must}
The measures \(\mu^\vr\) are stochastically nondecreasing in \(\vr\).
\end{cor}
\begin{proof}
Since $\vr$ and $\te$ are strictly increasing functions of each other, it 
is equivalent to show
monotonicity of \(\mu^\te\). This follows if we can show \(0\le\frac{\di}{\di\te}\Ev^\te(\vp(\om))\)
 for an arbitrary bounded nondecreasing function \(\vp\). Lemma \ref{lm:dercov} transforms this derivative into the covariance of \(\vp(\om)\) and \(\om\), which is non-negative due to \(\vp\) being nondecreasing.  \end{proof}

Monotonicity of \(\wih\mu^\vr\) requires somewhat more of a convexity argument.

\begin{pr}
The family of measures \({\wih\mu}^\vr\), defined in \eqref{eq:muhat}, is
stochastically nondecreasing in \(\vr\).
\label{pr:muhatpr}\end{pr}
\begin{proof}
Start by rewriting the definition:
\[
\ba
{\wih\mu}^\vr(y)&=\frac{\Ev^\vr\bigl([\om-\vr]\cdot{\bf1}\{\om>y\}\bigr)}{\Vv^\vr(\om)}=\frac{\Cov^\vr(\om,\,{\bf1}\{\om>y\})}{\Cov^\vr(\om,\,\om)}\\
&=\left.\frac{\frac{\di}{\di\te}\Pv^\te\{\om>y\}}{\frac{\di}{\di\te}\vr(\te)}\right|_{\te=\te(\vr)}=\frac{\di}{\di\vr}\Pv^\vr\{\om>y\}.
\ea
\]
Let us denote the \({\wih\mu}^\vr\)-expectation 
by \({\wih\Ev}^\vr\).  Fix a bounded nondecreasing function \(\vp\). We need to show
\[
0\le\frac{\di}{\di\vr}{\wih\Ev}^\vr\vp(\om).
\]

We compute a different expression for this derivative.  Passing the
derivative through the sum in the third equality below is 
justified because the series involved are dominated by   
certain geometric series, uniformly over $\te$  in small 
open neighborhoods.  This follows from the definitions of
 \(\un\te\) and \(\bar\te\) and the assumption \(\un\te<\te(\vr)<\bar\te\).
\[
\ba
{\wih\Ev}^\vr\vp(\om)&=\sum_{y=\omin}^{\omax}\vp(y)\cdot\frac{\di}{\di\vr}\Pv^\vr\{\om>y\}\\
&=\sum_{y=\omin}^{\omax}\vp(y)\cdot\frac{\di}{\di\vr}[\Pv^\vr\{\om>y\}-{\bf1}\{0\ge y\}]\\
&=\frac{\di}{\di\vr}\sum_{y=\omin}^{\omax}\vp(y)\cdot[\Pv^\vr\{\om>y\}-{\bf1}\{0\ge y\}]\\
&=\frac{\di}{\di\vr}\Ev^\vr\sum_{y=\omin}^{\omax}\vp(y)\cdot[{\bf1}\{\om>y\}-{\bf1}\{0\ge y\}]\\
&=\frac{\di}{\di\vr}\Ev^\vr\sum_{y=\omin}^{\omax}\vp(y)\cdot[{\bf1}\{\om>y>0\}-{\bf1}\{0\ge y\ge\om\}]\\
&=\frac{\di}{\di\vr}\Ev^\vr\Bigl[\sum_{y=1}^{\om-1}\vp(y)-\sum_{y=\om}^0\vp(y)\Bigr]=\frac{\di}{\di\vr}\Ev^\vr\Phi(\om).
\ea
\]
Above we introduced 
 the function \[\Phi(x)=\sum\limits_{y=1}^{x-1}\vp(y)-
\sum\limits_{y=x}^0\vp(y), \]
 with the convention that empty sums are zero. To conclude the proof, notice that \(\Phi(x+1)-\Phi(x)=\vp(x)\). 
Thus a nondecreasing function \(\vp\) determines a (non-strictly) 
convex function \(\Phi\) with \(\Phi(1)=0\), and vice-versa. 
Hence Theorem \ref{tm:convtm1} establishes that
\[
\frac{\di}{\di\vr}{\wih\Ev}^\vr\vp(\om)=\frac{\di^2}{\di\vr^2}\Ev^\vr\Phi(\om)\ge0.
\qedhere
\]
\end{proof}

\section{Regularity properties of the hydrodynamic flux function}\label{sc:afx}

For the zero range process defined among the examples  in Section \ref{sc:genex}, the hydrodynamic
(macroscopic)  
flux function \(\mathcal H\,:\,\Rb^+\to\Rb^+\) of \eqref{eq:fluxdef} is given by 
\[
\mathcal H(\vr) =\Ev^\vr f(\om).
\]
The results of 
Section \ref{sc:convtotpos} for \(f\)  now read as follows:
\begin{pr}\label{pr:zrcvx}
If the jump rate \(f\) of the zero range process is convex (or concave), then the  
 flux \(\mathcal H\) is also convex (or concave, respectively). Moreover, in this case \(\mathcal H''(\vr)>0\) 
 (or $\Hc''(\vr)<0$, respectively) for all $\vr>0$  if and only if \(f\) is not a linear function.
\end{pr}

Parts of this proposition were proved with  coupling methods in \cite{fluct}.

Next we show in the general case  that   $\Hc(\vr)$  is well defined, and
 is infinitely differentiable. (We use  third derivatives in
the proof of Theorem \ref{tm:main}.) 
The function $\Hc(\vr)$ is, in general, the expected net growth rate w.r.t.\ \(\mu^\vr\) as defined in \eqref{eq:fluxdef}.
We show that the series making up this expectation is finite, even after differentiating its terms. This will then lead to
smoothness of \(\Hc(\vr)\).
\begin{lm}\label{lm:hreg}
Let \(g(y,\,z)\ge0\) be any function on \(\Zb\times\Zb\), bounded by a polynomial in \(|y|\) and \(|z|\). Then for any \(\un\te<\te<\bar\te\),
\[
\Ev^\te\bigl[(p(\om_0,\,\om_1)+q(\om_0,\,\om_1))g(\om_0,\,\om_1)\bigr]<\infty.
\]
\end{lm}
\begin{proof}
We deal with the first part that contains \(p\), the one with \(q\) can be treated analogously. The sum we are looking at is\[
\sum_{y=\omin+1}^{\omax}\sum_{z=\omin}^{\omax-1}p(y,\,z)\cdot g(y,\,z)\cdot\frac{e^{\te(y+z)}}{f(y)!\cdot f(z)!}\cdot\frac{1}{Z(\te)^2}.
\]
These sums are certainly convergent if $\omin$ and $\omax$ are both finite. When this is not the case we split both summations at zero, and convergence is established on the four quadrants of the plane. We use \eqref{eq:mon} and the corollary
\[
p(y,\,z)=p(z+1,\,y-1)\cdot\frac{f(y)}{f(z+1)}\qquad\text{for }\omin<y\le\omax\text{ and }\omin\le z<\omax
\]
of \eqref{eq:symm}, and we consider empty sums to be zero.
\begin{itemize}
\item \(y>0,\,z>0\): In this case
\[
p(y,\,z)\le p(y,\,0)=p(1,\,y-1)\cdot\frac{f(y)}{f(1)}\le p(1,\,0)\cdot\frac{f(y)}{f(1)},
\]
and the corresponding part of the summation is bounded by
\[
\frac{p(1,\,0)}{f(1)}\cdot\sum_{y=1}^{\omax}\sum_{z=1}^{\omax-1}g(y,\,z)\cdot\frac{e^{\te(y+z)}}{f(y-1)!\cdot f(z)!}\cdot\frac{1}{Z(\te)^2}.
\]
\item \(y\le0,\,z>0\): In this case
\[
p(y,\,z)\le p(1,\,0),
\]
and the corresponding part of the summation is bounded by
\[
p(1,\,0)\cdot\sum_{y=\omin+1}^0\sum_{z=1}^{\omax-1}g(y,\,z)\cdot\frac{e^{\te(y+z)}}{f(y)!\cdot f(z)!}\cdot\frac{1}{Z(\te)^2}.
\]
\item \(y\le0,\,z\le0\): In this case
\[
p(y,\,z)\le p(1,\,z)=p(z+1,\,0)\cdot\frac{f(1)}{f(z+1)}\le p(1,\,0)\cdot\frac{f(1)}{f(z+1)},
\]
and the corresponding part of the summation is bounded by
\[
p(1,\,0)f(1)\cdot\sum_{y=\omin+1}^0\sum_{z=\omin}^0g(y,\,z)\cdot\frac{e^{\te(y+z)}}{f(y)!\cdot f(z+1)!}\cdot\frac{1}{Z(\te)^2}.
\]
\item \(y>0,\,z\le0\): In this case
\[
p(y,\,z)=p(z+1,\,y-1)\cdot\frac{f(y)}{f(z+1)}\le p(1,\,0)\cdot\frac{f(y)}{f(z+1)},
\]
and the corresponding part of the summation is bounded by
\[
p(1,\,0)\cdot\sum_{y=1}^{\omax}\sum_{z=\omin}^0g(y,\,z)\cdot\frac{e^{\te(y+z)}}{f(y-1)!\cdot f(z+1)!}\cdot\frac{1}{Z(\te)^2}.
\]
\end{itemize}
Convergence of each of these bounds for \(\un\te<\te<\bar\te\) is established e.g.\ by the ratio test.
\end{proof}
Notice that a similar argument gives finite higher moments of the rates when \(\log(f)\) is at most linear in both directions on \(\Zb\).
\begin{cor}
  $\Hc(\vr)$ is infinitely   differentiable at all \(\vr\) $\in$ (\(\omin,\,\omax)\).
\end{cor}\begin{proof}
By the previous lemma the series 
\[
F(\te):\,=\Hc(\vr(\te))=\frac1{Z(\te)^2}\cdot\sum_{y,\,z=\omin}^{\omax}(p(y,\,z)-q(y,\,z))\frac{e^{\te(y+z)}}{f(y)!\cdot f(z)!},
\]
is convergent and infinitely   differentiable. Since
$
\Hc(\vr)=F(\te(\vr))
$
and \(\vr\mapsto\te(\vr)\) is infinitely   differentiable as well, the claim follows.\end{proof}

\section*{Acknowledgment}
We thank an anonymous referee for careful reading of the manuscript and valuable suggestions.

\bibliography{refsmarton}
\bibliographystyle{plain}

\end{document}